# ASYMPTOTIC PROPERTIES OF THE MAXIMUM LIKELIHOOD ESTIMATOR IN AUTOREGRESSIVE MODELS WITH MARKOV REGIME

By Randal Douc[1], Éric Moulines[1] and Tobias Rydén[2]

*Ecole Nationale Supérieure des Télécommunications, Ecole Nationale Supérieure des Télécommunications and Lund University*

An autoregressive process with Markov regime is an autoregressive process for which the regression function at each time point is given by a nonobservable Markov chain. In this paper we consider the asymptotic properties of the maximum likelihood estimator in a possibly nonstationary process of this kind for which the hidden state space is compact but not necessarily finite. Consistency and asymptotic normality are shown to follow from uniform exponential forgetting of the initial distribution for the hidden Markov chain conditional on the observations.

**1. Introduction.** An autoregressive process with Markov regime, or Markov-switching autoregression, is a bivariate process $\{(X_k, Y_k)\}$, where $\{X_k\}$ is a Markov chain on a state space $\mathcal{X}$ and, conditional on $\{X_k\}$, $\{Y_k\}$ is an inhomogeneous $s$-order Markov chain on a state space $\mathcal{Y}$ such that the conditional distribution of $Y_n$ only depends on $X_n$ and lagged $Y$'s. The process $\{X_k\}$, usually referred to as the *regime*, is not observable and inference has to be carried out in terms of the observable process $\{Y_k\}$. In general we can write a model of this kind as

$$Y_n = f_\theta(\overline{\mathbf{Y}}_{n-1}, X_n; e_n),$$

where $\{e_k\}$ is an independent and identically distributed sequence of random variables that we denote the innovation process (the $e$'s are not the

Received May 2001; revised September 2003.
[1]Supported by the EU TMR network Statistical and Computational Methods for the Analysis of Spatial Data.
[2]Supported by a grant from the Swedish Research Council for Engineering Sciences.
*AMS 2000 subject classifications.* Primary 62M09; secondary 62F12.
*Key words and phrases.* Asymptotic normality, autoregressive process, consistency, geometric ergodicity, hidden Markov model, identifiability, maximum likelihood, switching autoregression.

innovation process in Wold's sense, however), $\overline{\mathbf{Y}}_k \triangleq (Y_k, Y_{k-1}, \ldots, Y_{k-s+1})$ and $\{f_\theta\}$ is a family of functions indexed by a finite-dimensional parameter $\theta$. Of particular interest are the linear autoregressive models for which

$$f_\theta(\overline{\mathbf{Y}}_{n-1}, X_n; e_n) = \sum_{i=1}^{s} a_i(X_n; \theta) Y_{n-i} + e_n.$$

These models were initially proposed by Hamilton (1989) in econometric theory; the number of states of the Markov chain is in this context most often assumed to be finite, each state being associated with a given state of the economy [see Krolzig (1997), Kim and Nelson (1999) and references therein]. Linear autoregressive processes with Markov regime are also widely used in several electrical engineering areas including tracking of maneuvering targets [Bar-Shalom and Li (1993)], failure detection [Tugnait (1982)] and stochastic adaptive control [Doucet, Logothetis and Krishnamurthy (2000)]; in such cases the hidden state is most often assumed to be continuous. Nonlinear switching autoregressive models have recently been proposed in quantitative finance to model volatility of log-returns of international equity markets [see, e.g., Susmel (2000) and Chib, Nardari and Shephard (2002)]. A simple example of such a model (referred to as SWARCH for switching ARCH) is

$$Y_n = f_\theta(\overline{\mathbf{Y}}_{n-1}, X_n) e_n,$$

where once again $\{X_k\}$ is either a discrete or a continuous Markov chain. Another important subclass of autoregressive models with Markov regime are the hidden Markov models (HMMs), for which the conditional distribution of $Y_n$ does not depend on lagged $Y$'s but only on $X_n$. HMMs are used in many different areas, including speech recognition [Juang and Rabiner (1991)], neurophysiology [Fredkin and Rice (1987)], biology [Churchill (1989)], econometrics [Chib, Nardari and Shephard (2002)] and time series analysis [de Jong and Shephard (1995) and Chan and Ledolter (1995)]. See also the monograph by MacDonald and Zucchini (1997) and references therein.

Most works on maximum likelihood estimation in such models have focused on numerical methods suitable for approximating the maximum likelihood estimator (MLE). In sharp contrast, statistical issues regarding asymptotic properties of the MLE for autoregressive models with Markov regime have been largely ignored until recently. Baum and Petrie (1966) proved consistency and asymptotic normality of the MLE for HMMs in the particular case where both the observed and the latent variables take values is finite spaces. These results have recently been extended in a series of papers by Leroux (1992), Bickel and Ritov (1996), Bickel, Ritov and Rydén (1998) (henceforth referred to as BRR), Jensen and Petersen (1999) (henceforth referred to as JP) and Bakry, Milhaud and Vandekerkhove (1997). BRR followed the approach taken by Baum and Petrie (1966) and generalized their results to the case where the hidden Markov chain $\{X_k\}$ takes a finite number of values, but the observations belong to a general space. JP extended these results to HMMs with the regime taking values in a compact space, proving asymptotic normality of the MLE and a local consistency theorem.

Around the same time, Le Gland and Mevel (2000) [see also Mevel (1997)] independently developed a different technique to prove consistency and asymptotic normality of the MLE for HMMs with finite hidden state space. Their work was later extended to HMMs with nonfinite hidden state space by Douc and Matias (2001). This approach is based on the observation that the log likelihood can be expressed as an additive function of an extended Markov chain. These techniques, which are well adapted to study recursive estimators (that are updated for each novel observation), typically require stronger assumptions than the methods developed in BRR and JP.

None of the theoretical contributions mentioned so far allows for autoregression, but are concerned with HMMs alone. For autoregressive processes with Markov regime, the only theoretical result available up till now is consistency of the MLE when the regime takes values in a finite set [Krishnamurthy and Rydén (1998) and Francq and Roussignol (1998)]. In the present paper we examine asymptotic properties of the MLE when the hidden Markov chain takes values in a compact space, and we do allow for autoregression in the observable process. Our results include consistency and asymptotic normality of the MLE under standard regularity assumptions (Theorems 1 and 4) and consistency of the observed information as an estimator of the Fisher information (Theorem 3 with $\theta_n^*$ being the MLE). These results generalize what is obtained in the above-mentioned papers to a larger class of models, and we obtain them through a unified approach. We also point out that the convergence theorem for the MLE is global, as opposed to the local theorem of JP. Moreover, the nonstationary setting is treated in Section 7.

The likelihood that we will work with is the conditional likelihood given initial observations $\overline{\mathbf{Y}}_0 = (Y_0, \ldots, Y_{-s+1})$ and the initial (but unobserved) state $X_0$. Conditioning on initial observations in time series models goes back at least to Mann and Wald (1943). In our case we, in addition, also condition the likelihood on the unobserved initial state. The reason for doing so is that the stationary distribution of $\{(X_k, Y_k)\}$, and hence the true likelihood, is typically infeasible to compute. Thus $n$, denoting the number of factors in the likelihood—the "nominal" sample size—is $s$ less than the actual sample size. Using $p$ as a generic symbol for densities we can express the conditional log likelihood as

$$\log p_\theta(y_1, \ldots, y_n | \bar{\mathbf{y}}_0, x_0)$$
$$(1) \qquad = \log \int \cdots \int p_\theta(x_1, \ldots, x_n, y_1, \ldots, y_n | \bar{\mathbf{y}}_0, x_0) \mu(dx_1) \cdots \mu(dx_n)$$

$$= \log \int \cdots \int \prod_{k=1}^{n} q_\theta(x_{k-1}, x_k) \prod_{k=1}^{n} g_\theta(y_k|\bar{\mathbf{y}}_{k-1}, x_k) \mu(dx_1) \cdots \mu(dx_n),$$

where $\mu$ and $q_\theta(\cdot, \cdot)$ are a reference measure and the transition density for the hidden chain, respectively, and $g_\theta(y_k|\bar{\mathbf{y}}_{k-1}, x_k)$ is the conditional density of $y_k$ given $\bar{\mathbf{y}}_{k-1}$ and $x_k$. In the particular case when $\{X_k\}$ is finite-valued, taking values in $\{1, 2, \ldots, d\}$ say, this log likelihood can be expressed as

$$(2) \qquad \log p_\theta(y_1, \ldots, y_n | \bar{\mathbf{y}}_0, x_0) = \log \mathbf{1}_{x_0}^T \left( \prod_{k=1}^{n} Q_\theta G_\theta(y_k | \bar{\mathbf{y}}_{k-1}) \right) \mathbf{1},$$

where $Q_\theta = \{q_\theta(i, j)\}$ is the transition probability matrix of the Markov chain $\{X_k\}$, $G_\theta(y|\bar{\mathbf{y}}) = \mathrm{diag}(g_\theta(y|\bar{\mathbf{y}}, i))$, $\mathbf{1}_{x_0}$ is the $x_0$th unit vector of length $d$, that is, a $d \times 1$ vector in which all elements are zero except for element $x_0$ which is unity, and $\mathbf{1}$ is a $d \times 1$ vector of all ones. It is clear that (2) is essentially a product of matrices and is hence easily evaluated. It can be maximized over $\theta$ using standard numerical optimization procedures or using the EM algorithm [see, e.g., Hamilton (1990)]. However, one should be aware that the log likelihood is typically multi-modal and either approach may converge to a local maximum. When $\{X_k\}$ is continuous, evaluation of the log likelihood (1) requires an integration over an $n$-dimensional space. This task is insurmountable for typical values of $n$, and approximation methods are required. Two classes of such methods, particle filters and Monte Carlo EM algorithms, as well as a numerical example using the latter, are briefly discussed in Section 8.

An obvious variant to our approach is to replace the condition of a fixed $x_0$ by assuming a fixed distribution for $x_0$. Such an assumption does not change any of our results and no more than notational changes are needed in the proofs. A further natural variant is to maximize (1) w.r.t. $\theta$ *and* the unknown $x_0$. We have not included this approach in the present paper, primarily because score function analysis would require assumptions on how the maximizing $x_0$ varies with $\theta$, assumptions that would be difficult to verify in practice. We do remark, however, that in a particular but important case, assuming a fixed $x_0$ is no less general than is maximization over $x_0$. Suppose that the regime $\{X_k\}$ is finite-valued and that all elements $q_{ij}$ of the transition probability matrix $Q$ may be chosen independently. The parameter vector $\theta$ may then be written $\theta = ((q_{ij}), \psi)$. We also assume that $\psi$ can be further decomposed as $\psi = (\alpha, (\beta_i))$, and that the functions $g$ are such that $g_\theta(y_k|\bar{\mathbf{y}}_{k-1}, x_k) = h(y_k|\bar{\mathbf{y}}_{k-1}; \alpha, \beta_{x_k})$ for some family of densities $h$. In other words, all $g$'s belong to a single parametric class of densities, $\alpha$ is a parameter common to all regimes and the $\beta_i$'s are the regime specific parameters. For example, in the linear regression case $\beta_i$ may be the regime specific regression coefficients while $\alpha$ may be a common innovation variance, $\alpha = \mathbb{E}e_n^2$. With this general structure it is clear that if $x_0$ is a fixed initial state, for any model with a different initial state we can find an equivalent model with initial state $x_0$ by simply renumbering the states and then reordering the $q_{ij}$'s and $\beta_i$'s accordingly. Therefore, whenever $\theta$ is structured as above, assuming a fixed $x_0$ is no less general than is maximization over $x_0$.

As mentioned above, from a practical point of view the novelty of the present paper is that we extend the analysis of MLE asymptotics to wider class of models using a unified approach. From a theoretical point of view the novelty is, foremost, the geometrically decaying bound on the mixing rate of the conditional chain, $X|Y$, given in Corollary 1 and (20). This bound parallels results of BRR (page 1622) and JP (page 521), but in contrast to those results our bound does not depend on the $Y$'s being conditioned upon; it is deterministic. Assumption (A1)(a) below, implying that the hidden chain is uniformly geometrically ergodic, and more specifically that the whole state space is 1-small [see the comment after (A1)(a)], is crucial to this property; if $\{X_k\}$ is $m$-small with $m > 1$ one can prove an analogous mixing rate bound using similar ideas, but the bound will then depend on the $Y$'s. The deterministic nature of the bound is vital to our proofs that the conditional score given the "infinite past" [$\Delta_{k,\infty}(\theta^*)$ in Section 6.1] and the conditional Hessian given the "infinite past" (cf. Propositions 4 and 5) have finite second and first moments, respectively. The reason is that when the model contains autoregression, the conditional distribution of $\{Y_k\}$ given $\{X_k\}$ is governed by an inhomogeneous autoregression rather than by independence; hence, in the proof of Lemma 10, for example, we cannot condition on the regime $\{X_k\}$ and exploit conditional independence in order to turn a random mixing bound into a deterministic one as was done in BRR (e.g., page 1625) and JP (e.g., page 525). We plan to look into this more general case, but it lies outside the scope of the present paper. Another feature of the present paper is that by refining the arguments of BRR and JP we obtain almost sure convergence rather than convergence is probability in Theorem 3.

The paper is organized as follows. Main assumptions are given and commented in Section 2, together with common notation. Then in Section 3 we show that the regime $\{X_k\}$, given the observations, is a nonhomogeneous Markov chain whose transition kernels may be minorized using a fixed and common minorizing constant. This leads to a deterministic bound for its mixing rate. In Section 4, consistency of the MLE is considered under the additional assumption that $\{Y_k\}$ is strict sense stationary; extensions to nonstationary processes through coupling are carried out in Section 7. Conditions upon which the parameters are identifiable are given in Section 5. Asymptotic normality of the estimator is studied in Section 6. The proof is based on a central limit theorem and a locally uniform law of large numbers for the conditional expectation of appropriately defined statistics. More specifically, these statistics are additive and quadratic functionals of the complete data. Section 8 contains a discussion of numerical methods for state space models and a numerical example. Finally, the Appendix contains proofs not given in the main text.

**2. Notation and assumptions.** We assume that the Markov chain $\{X_k\}_{k=0}^\infty$ is homogeneous and lies in a separable and compact set $\mathcal{X}$, equipped with a metrizable topology and the associated Borel $\sigma$-field $\mathcal{B}(\mathcal{X})$. We let $Q_\theta(x, A)$, $x \in \mathcal{X}$, $A \in \mathcal{B}(\mathcal{X})$, be the transition kernel of the chain; the parameter $\theta$ which indexes the family of transition kernels as well as the regression functions for the $Y$'s, see below, is the parameter that we want to estimate. Next we assume that each measure $Q_\theta(x, \cdot)$ has a density $q_\theta(x, \cdot)$ with respect to a common *finite* dominating measure $\mu$ on $\mathcal{X}$. That is, for all $\theta$ and $x \in \mathcal{X}$, $Q_\theta(x, \cdot) \ll \mu$. For the sake of simplicity, it is assumed that $\mu(\mathcal{X}) = 1$; this assumption hints at applications where $\mathcal{X}$ is a totally bounded space.

We also assume that the observable sequence $\{Y_k\}_{k=-s+1}^\infty$ takes values in a set $\mathcal{Y}$ that is separable and metrizable by a complete metric. Furthermore, for each $n \geq 1$ and given $\{Y_k\}_{k=n-s}^{n-1}$ and $X_n$, $Y_n$ is conditionally independent of $\{Y_k\}_{k=-s+1}^{n-s-1}$ and $\{X_k\}_{k=0}^{n-1}$. We also assume that for each $X_n$, $\overline{Y}_{n-1}$ and $\theta$, this conditional law has a density $g_\theta(y|\overline{Y}_{n-1}, X_n)$ with respect to some fixed $\sigma$-finite measure $\nu$ on the Borel $\sigma$-field $\mathcal{B}(\mathcal{Y})$.

The parameter $\theta$ belongs to $\Theta$, a compact subset of $\mathbb{R}^p$. The true parameter value will be denoted by $\theta^*$, and when proving asymptotic normality of the MLE we assume that $\theta^*$ lies in the interior of $\Theta$. Given the observations $Y_{-s+1}, \ldots, Y_n$ of the process $\{Y_k\}$, we wish to estimate $\theta^*$ by the maximum likelihood method.

The sequence $\{Z_k\}_{k=0}^\infty \triangleq \{(X_k, \overline{Y}_k)\}_{k=0}^\infty$ is a Markov chain on $\mathcal{Z} \triangleq \mathcal{X} \times \mathcal{Y}^s$ with transition kernel $\Pi_\theta$ given by, for any bounded measurable function $f$ on $\mathcal{Z}$,

$$\Pi_\theta f(x, y_s, y_{s-1}, \ldots, y_1)$$
$$= \int_{\mathcal{X} \times \mathcal{Y}} f(x', y', y_s, \ldots, y_2) q_\theta(x, x') g_\theta(y'|y_s, \ldots, y_1, x') \mu(dx') \nu(dy').$$

We use in the sequel the canonical version of this Markov chain and put $\bar{\nu} \triangleq \nu^{\otimes s}$. For a probability measure $\zeta$ on $\mathcal{Z}$ we let $\mathbb{P}_{\theta,\zeta}$ be the law of $\{Z_n\}$ when the initial distribution is $\zeta$; that is, $Z_0 \sim \zeta$. Furthermore, $\mathbb{E}_{\theta,\zeta}$ is the associated expectation. Many conditional probabilities and expectations in this paper do not depend on the initial distribution, and we stress this by then dropping the initial probability measure from the notation, so that $\mathbb{P}_{\theta,\zeta}$ is replaced by $\mathbb{P}_\theta$, and so on.

Throughout this paper we will assume that the transition kernel $\Pi_\theta$ has a unique invariant distribution $\pi_\theta$; this assumption is further commented on below. For a stationary process we write $\overline{\mathbb{P}}_\theta$ and $\overline{\mathbb{E}}_\theta$ for $\mathbb{P}_{\theta,\pi_\theta}$ and $\mathbb{E}_{\theta,\pi_\theta}$, respectively. We can and will extend such a stationary process $\{Z_k\}_{k=0}^{\infty}$ to a stationary Markov chain $\{Z_k\}_{k=-\infty}^{\infty}$ with doubly infinite time and the same transition kernel.

For $i \leq j$, put $\mathbf{Y}_i^j \triangleq (Y_i, Y_{i+1}, \ldots, Y_j)$ and $\overline{\mathbf{Y}}_i^j \triangleq (\overline{\mathbf{Y}}_i, \overline{\mathbf{Y}}_{i+1}, \ldots, \overline{\mathbf{Y}}_j)$, respectively. Similar notation will be used for other quantities. For any measurable function $f$ on $(\mathcal{X}, \mathcal{B}(\mathcal{X}), \mu)$, $\operatorname{ess\,sup} f \triangleq \inf\{M \geq 0 : \mu(\{M \leq |f|\}) = 0\}$ and, if $f$ is nonnegative, $\operatorname{ess\,inf} f \triangleq \sup\{M \geq 0 : \mu(\{M \geq f\}) = 0\}$ (with obvious conventions if those sets are empty). For the sake of simplicity, instead of writing $\operatorname{ess\,sup} f$ or $\operatorname{ess\,inf} f$, we use the notation $\sup f$ or $\inf f$. For any two probability measures $\mu_1$ and $\mu_2$ we define the total variation distance $\|\mu_1 - \mu_2\|_{\mathrm{TV}} = \sup_A |\mu_1(A) - \mu_2(A)|$ and we also recall the identities $\sup_{|f|\leq 1} |\mu_1(f) - \mu_2(f)| = 2\|\mu_1 - \mu_2\|_{\mathrm{TV}}$ and $\sup_{0 \leq f \leq 1} |\mu_1(f) - \mu_2(f)| = \|\mu_1 - \mu_2\|_{\mathrm{TV}}$. For any matrix or vector $A$, $\|A\| = \sum |A_{ij}|$. Finally, we will use the letter $p$ to denote densities w.r.t. the probability measure on $\mathcal{B}(\mathcal{X} \times \mathcal{Y})^{\otimes \mathbb{Z}}$ whose finite-dimensional distributions are given by $(\mu \otimes \nu)^{\otimes r}$ for all $r \geq 1$.

We now list our basic assumptions.

(A1) (a) $0 < \sigma_- \triangleq \inf_{\theta \in \Theta} \inf_{x,x' \in \mathcal{X}} q_\theta(x, x')$ and $\sigma_+ \triangleq \sup_{\theta \in \Theta} \sup_{x,x' \in \mathcal{X}} q_\theta(x, x') < \infty$.
 (b) For all $y' \in \mathcal{Y}$ and $\bar{\mathbf{y}} \in \mathcal{Y}^s$, $0 < \inf_{\theta \in \Theta} \int_{\mathcal{X}} g_\theta(y'|\bar{\mathbf{y}}, x) \mu(dx)$ and $\sup_{\theta \in \Theta} \int_{\mathcal{X}} g_\theta(y'|\bar{\mathbf{y}}, x) \mu(dx) < \infty$.

Assumption (A1)(a) implies that for all $x \in \mathcal{X}$, $Q(x, A) \geq \sigma_- \mu(A)$ where $\mu$ is a probability measure, that is, the state space $\mathcal{X}$ of the Markov chain $\{X_n\}$ is 1-small [Meyn and Tweedie (1993), page 106, with $m = 1$]. Thus, for all $\theta \in \Theta$, this chain has a unique invariant measure $\pi_\theta^X$ and is uniformly ergodic [Meyn and Tweedie (1993), Theorem 16.0.2(v)]. When the state space is finite, (A1)(a) is equivalent to saying that for all $x, x' \in \mathcal{X}$, $\inf_{\theta \in \Theta} q_\theta(x, x') > 0$.

(A2) For all $\theta \in \Theta$, the transition kernel $\Pi_\theta$ is positive Harris recurrent and aperiodic with invariant distribution $\pi_\theta$.

That the chain is positive means, essentially, that it is irreducible and has an invariant distribution [Meyn and Tweedie (1993), page 230] and Harris recurrence means that any nonnull set will be infinitely visited by the chain irrespective of where it starts within the set [Meyn and Tweedie (1993), page 200]. This assumption is rather weak; results on ergodicity for autoregressive processes with Markov regime can be found in, for example, Francq and Zakoian (2001), Holst, Lindgren, Holst and Thuvesholmen [(1994), page 495] and Yao and Attali (2000). It implies that for any initial measure $\lambda$ [see Meyn and Tweedie (1993), Theorem 13.3.3],

$$\lim_{n \to \infty} \|\lambda \Pi_\theta^n - \pi_\theta\|_{\mathrm{TV}} = 0, \tag{3}$$

so that the tail $\sigma$-field of $\{Z_k\}$ is trivial [Lindvall (1992), Theorem III.21.12]. Its invariant $\sigma$-field, which is no larger, is thus also trivial and hence $\{Z_k\}$ is ergodic in the measure-theoretic sense of the word.

For the developments that follow, an additional assumption is needed.

(A3) $b_+ \triangleq \sup_\theta \sup_{\bar{\mathbf{y}}_0, y_1, x} g_\theta(y_1|\bar{\mathbf{y}}_0, x) < \infty$ and $\overline{\mathbb{E}}_{\theta^*}(|\log b_-(\overline{\mathbf{Y}}_0, Y_1)|) < \infty$, where $b_-(\bar{\mathbf{y}}_0, y_1) \triangleq \inf_\theta \int_{\mathcal{X}} g_\theta(y_1|\bar{\mathbf{y}}_0, x) \mu(dx)$.

REMARK 1. In the sequel we consider conditional expectations of random variables w.r.t. the $\sigma$-algebra generated by $(\mathbf{X}_m^n, \mathbf{Y}_m^n)$ for some $m \leq n$. Such expectations are defined up to a $\mathbb{P}_{\theta,\zeta}$-null set. For the derivations that follow, we need to specify a version of these conditional expectations. Since $\mathbb{P}_{\theta,\zeta}$ is defined by the initial distribution $\zeta$ and the transition kernel $\Pi_\theta$, it is always possible to express these conditional expectations in terms of these quantities and we always implicitly choose this version of the conditional expectations.

**3. Uniform forgetting of the conditional hidden Markov chain.** By the conditional hidden Markov chain we mean the process $\{X_k\}$ given a sequence of $Y$'s. It will turn out that this process is a Markov chain, although nonhomogeneous, but still having a favorable mixing rate. Bounds on this mixing rate will be instrumental in the forthcoming development.

LEMMA 1. *Assume* (A1). *Let* $m, n \in \mathbb{Z}$ *with* $m \leq n$ *and* $\theta \in \Theta$. *Under* $\overline{\mathbb{P}}_\theta$, *conditionally on* $\overline{\mathbf{Y}}_m^n$, $\{X_k\}_{k \geq m}$ *is an inhomogeneous Markov chain, and for all* $k > m$ *there exists a function* $\mu_k(\mathbf{y}_{k-s}^n, A)$ *such that*:

(i) *for any* $A \in \mathcal{B}(\mathcal{X})$, $\mathbf{y}_{k-s}^n \mapsto \mu_k(\mathbf{y}_{k-s}^n, A)$ *is a Borel function*;
(ii) *for any* $\mathbf{y}_{k-s}^n$, $\mu_k(\mathbf{y}_{k-s}^n, \cdot)$ *is a probability measure on* $\mathcal{B}(\mathcal{X})$. *In addition, for all* $\mathbf{y}_{k-s}^n$ *it holds that* $\mu_k(\mathbf{y}_{k-s}^n, \cdot) \ll \mu$ *and for all* $\overline{\mathbf{Y}}_m^n$,

$$\inf_{x \in \mathcal{X}} \overline{\mathbb{P}}_\theta(X_k \in A | X_{k-1} = x, \overline{\mathbf{Y}}_m^n) \geq \frac{\sigma_-}{\sigma_+} \mu_k(\mathbf{Y}_{k-s}^n, A).$$

REMARK 2. Contrary to JP, this minorization condition involves a constant $\sigma_-/\sigma_+$ which does not depend on the values of $\{Y_k\}$. On the other hand, the minorizing measure $\mu_k(\mathbf{y}_{k-s}^n, \cdot)$ does depend on $\mathbf{y}_{k-s}^n$ whereas the minorizing measure is fixed in JP. Hence no assumption on the conditional density of $Y_k$ given past observations and hidden state variables is needed, whereas JP assumed a moment condition, in the special case of HMMs, for the ratio $\sup_\theta \sup_{x,x'} g_\theta(y|x)/g_\theta(y|x')$. An explicit expression for $\mu_k(\mathbf{y}_{k-s}^n, \cdot)$ is not needed.

PROOF OF LEMMA 1. The proof is adapted from Del Moral and Guionnet (2001) [see also Del Moral and Miclo (2000)]. The Markov property implies that, for $m < k \leq n$,

$$\overline{\mathbb{P}}_\theta(X_k \in A | \mathbf{X}_m^{k-1}, \overline{\mathbf{Y}}_m^n) = \overline{\mathbb{P}}_\theta(X_k \in A | X_{k-1}, \overline{\mathbf{Y}}_{k-1}^n).$$

For $k > n$ we have $\overline{\mathbb{P}}_\theta(X_k \in A | \mathbf{X}_m^{k-1}, \overline{\mathbf{Y}}_m^n) = Q_\theta(X_{k-1}, A)$. This shows that $\{X_k\}_{k \geq m}$ conditional on $\overline{\mathbf{Y}}_m^n$ is an inhomogeneous Markov chain. For $k \leq n$ it holds that

$$\overline{\mathbb{P}}_\theta(X_k \in A | X_{k-1}, \overline{\mathbf{Y}}_{k-1}^n)$$
$$= \int_A q_\theta(X_{k-1}, x) \bar{p}_\theta(\mathbf{Y}_k^n | X_k = x, \overline{\mathbf{Y}}_{k-1}) \mu(dx)$$
$$\times \left( \int_\mathcal{X} q_\theta(X_{k-1}, x) \bar{p}_\theta(\mathbf{Y}_k^n | X_k = x, \overline{\mathbf{Y}}_{k-1}) \mu(dx) \right)^{-1},$$

where

(4)
$$\bar{p}_\theta(\mathbf{Y}_k^n | X_k = x_k, \overline{\mathbf{Y}}_{k-1})$$
$$= \int \prod_{i=k+1}^n q_\theta(x_{i-1}, x_i) \prod_{i=k}^n g_\theta(Y_i | \overline{\mathbf{Y}}_{i-1}, x_i) \mu^{\otimes (n-k)}(d\mathbf{x}_{k+1}^n).$$

Since $\sigma_- \leq q_\theta(x, x') \leq \sigma_+$ it readily follows that

$$\overline{\mathbb{P}}_\theta(X_k \in A | X_{k-1}, \overline{\mathbf{Y}}_{k-1}^n) \geq \frac{\sigma_-}{\sigma_+} \mu_k(\mathbf{Y}_{k-s}^n, A)$$

with

$$\mu_k(\mathbf{Y}_k^n, A) \triangleq \int_A \bar{p}_\theta(\mathbf{Y}_k^n | X_k = x, \overline{\mathbf{Y}}_{k-1}) \mu(dx) \bigg/ \int_\mathcal{X} \bar{p}_\theta(\mathbf{Y}_k^n | X_k = x, \overline{\mathbf{Y}}_{k-1}) \mu(dx).$$

Note that

$$\int_\mathcal{X} \bar{p}_\theta(\mathbf{Y}_k^n | X_k = x_k, \overline{\mathbf{Y}}_{k-1}) \mu(dx_k)$$
$$= \int \prod_{i=k+1}^n q_\theta(x_{i-1}, x_i) \prod_{i=k}^n g_\theta(Y_i | \overline{\mathbf{Y}}_{i-1}, x_i) \mu^{\otimes (n-k+1)}(d\mathbf{x}_k^n)$$
$$\geq \sigma_-^{n-k} \prod_{i=k}^n \int g_\theta(Y_i | \overline{\mathbf{Y}}_{i-1}, x) \mu(dx)$$

is positive under (A1)(b). For $k > n$ we simply set $\mu_k(\mathbf{Y}_{k-s}^n, A) = \mu(A)$. □

The a posteriori chain thus also admits $\mathcal{X}$ as a 1-small set. It is worthwhile to note that, despite the chain being nonhomogeneous, the same minorizing constant can be used for all kernels, irrespective of the $Y$'s the chain is conditioned upon and of the parameter value. Using standard results for uniformly minorized Markov chains [see, e.g., Lindvall (1992), Sections III.9–11], we thus have the following result, which plays a key role in the sequel.

COROLLARY 1. *Assume* (A1). *Let $m, n \in \mathbb{Z}$ with $m \leq n$ and $\theta \in \Theta$. Then for all $k \geq m$, all probability measures $\mu_1$ and $\mu_2$ on $\mathcal{B}(\mathcal{X})$ and all $\overline{\mathbf{Y}}_m^n$,*

$$\left\| \int_{\mathcal{X}} \overline{\mathbb{P}}_\theta(X_k \in \cdot | X_m = x, \overline{\mathbf{Y}}_m^n) \mu_1(dx) - \int_{\mathcal{X}} \overline{\mathbb{P}}_\theta(X_k \in \cdot | X_m = x, \overline{\mathbf{Y}}_m^n) \mu_2(dx) \right\|_{\mathrm{TV}}$$
$$\leq \rho^{k-m},$$

*where $\rho \triangleq 1 - \sigma_-/\sigma_+$.*

Note that when $m$ is positive, $\overline{\mathbb{P}}_\theta(X_k \in \cdot | X_m = x, \overline{\mathbf{Y}}_m^n) = \mathbb{P}_\theta(X_k \in \cdot | X_m = x, \overline{\mathbf{Y}}_m^n)$ does not depend upon the initial distribution.

**4. Uniform convergence of the likelihood contrast function.** Given $x_0 \in \mathcal{X}$, notice that

$$(5) \quad p_\theta(\mathbf{Y}_1^n | \overline{\mathbf{Y}}_0, X_0 = x_0) = \int \prod_{k=1}^n q_\theta(x_{k-1}, x_k) g_\theta(Y_k | \overline{\mathbf{Y}}_{k-1}, x_k) \mu^{\otimes n}(d\mathbf{x}_1^n)$$

and define the *conditional* log likelihood function

$$(6) \quad l_n(\theta, x_0) \triangleq \log p_\theta(\mathbf{Y}_1^n | \overline{\mathbf{Y}}_0, X_0 = x_0) = \sum_{k=1}^n \log p_\theta(Y_k | \overline{\mathbf{Y}}_0^{k-1}, X_0 = x_0),$$

where $p_\theta(Y_k | \overline{\mathbf{Y}}_0^{k-1}, X_0 = x_0) = p_\theta(\mathbf{Y}_1^k | \overline{\mathbf{Y}}_0, X_0 = x_0) / p_\theta(\mathbf{Y}_1^{k-1} | \overline{\mathbf{Y}}_0, X_0 = x_0)$. With the notation introduced above, for $k \geq 1$,

$$p_\theta(Y_k | \overline{\mathbf{Y}}_0^{k-1}, X_0 = x_0)$$
$$(7) \qquad = \iint g_\theta(Y_k | \overline{\mathbf{Y}}_{k-1}, x_k) q_\theta(x_{k-1}, x_k)$$
$$\qquad\qquad \times \mathbb{P}_\theta(X_{k-1} \in dx_{k-1} | \overline{\mathbf{Y}}_0^{k-1}, X_0 = x_0) \mu(dx_k);$$

here $\mathbb{P}_\theta(X_{k-1} \in \cdot | \overline{\mathbf{Y}}_0^{k-1}, X_0 = x_0)$ is the filtering distribution of the unknown state $X_{k-1}$ given $\overline{\mathbf{Y}}_0^{k-1}$ and $X_0 = x_0$. Note that this distribution may be expressed as

$$(8) \quad \mathbb{P}_\theta(X_{k-1} \in \cdot | \overline{\mathbf{Y}}_0^{k-1}, X_0 = x) = \int \mathbb{P}_\theta(X_{k-1} \in \cdot | \overline{\mathbf{Y}}_0^{k-1}, X_0 = x_0) \delta_x(dx_0).$$

The discussion in the previous section hints that the influence of the initial point $X_0$ vanishes as $k \to \infty$.

The definition of the conditional log likelihood employed here differs from the one usually adopted for HMMs. Extending to AR models with Markov regime the definitions of BRR and JP for example, the log likelihood would be

$$l_n(\theta) \triangleq \sum_{k=1}^{n} \log \bar{p}_\theta(Y_k | \overline{\mathbf{Y}}_0^{k-1}), \tag{9}$$

where

$$\bar{p}_\theta(Y_k | \overline{\mathbf{Y}}_0^{k-1}) = \iint g_\theta(Y_k | \overline{\mathbf{Y}}_{k-1}, x_k) q_\theta(x_{k-1}, x_k) \\ \times \overline{\mathbb{P}}_\theta(X_{k-1} \in dx_{k-1} | \overline{\mathbf{Y}}_0^{k-1}) \mu(dx_k). \tag{10}$$

Here $\overline{\mathbb{P}}_\theta(X_{k-1} \in \cdot | \overline{\mathbf{Y}}_0^{k-1})$ is the filtering distribution of the unknown state $X_{k-1}$ given $\overline{\mathbf{Y}}_0^{k-1}$ under the stationary probability $\overline{\mathbb{P}}_\theta$. This filtering distribution may be expressed as

$$\overline{\mathbb{P}}_\theta(X_{k-1} \in \cdot | \overline{\mathbf{Y}}_0^{k-1}) = \int \mathbb{P}_\theta(X_{k-1} \in \cdot | \overline{\mathbf{Y}}_0^{k-1}, X_0 = x_0) \overline{\mathbb{P}}_\theta(X_0 \in dx_0 | \overline{\mathbf{Y}}_0^{k-1}) \tag{11}$$

and Corollary 1 shows that that the total variation distance between the filtering probabilities $\overline{\mathbb{P}}_\theta(X_{k-1} \in \cdot | \overline{\mathbf{Y}}_0^{k-1})$ and $\mathbb{P}_\theta(X_{k-1} \in \cdot | \overline{\mathbf{Y}}_0^{k-1}, X_0 = x_0)$ tends to zero exponentially fast as $k \to \infty$ uniformly w.r.t. $x_0$.

The definition of the log likelihood in (9) is useful for HMMs but less so for models with autoregression. Indeed, for many models $\bar{p}_\theta(Y_k | \mathbf{Y}_0^{k-1})$ cannot be expressed in closed form, basically because the smoothing probability $\overline{\mathbb{P}}_\theta(X_0 \in \cdot | \overline{\mathbf{Y}}_0^{k-1})$ depends upon the stationary distribution $\pi_\theta$ of the complete chain,

$$\overline{\mathbb{P}}_\theta(X_0 \in A | \overline{\mathbf{Y}}_0^{k-1}) \\ = \frac{\int_A \pi_\theta(dx_0 | \overline{\mathbf{Y}}_0) \int \prod_{i=1}^{k-1} q_\theta(x_{i-1}, x_i) g_\theta(Y_i | \overline{\mathbf{Y}}_{i-1}, x_i) \mu^{\otimes(k-1)}(d\mathbf{x}_1^{k-1})}{\int_\mathcal{X} \pi_\theta(dx_0 | \overline{\mathbf{Y}}_0) \int \prod_{i=1}^{k-1} q_\theta(x_{i-1}, x_i) g_\theta(Y_i | \overline{\mathbf{Y}}_{i-1}, x_i) \mu^{\otimes(k-1)}(d\mathbf{x}_1^{k-1})}.$$

In many models for which the stationary density is not available in closed form, the log likelihood (9) does not lead to a practical algorithm. This is our motivation for considering the conditional form (6) of the log likelihood function. Nevertheless, as we will see below, for any initial point $x_0$, $n^{-1}(l_n(\theta, x_0) - l_n(\theta))$ converges to zero uniformly w.r.t. to $\theta \in \Theta$ as a consequence of the uniform forgetting of the conditional Markov chain. Thus, by the continuity of the arg max functional, $\hat{\theta}_{n,x_0}$, the maximum of $l_n(\theta, x_0)$, and $\hat{\theta}_n$, the maximum of $l_n(\theta)$, are asymptotically equivalent.

REMARK 3. For $\xi$ an arbitrary probability measure on $\mathcal{B}(\mathcal{X})$ it is possible to consider
$$p_{\theta,\xi}(\mathbf{Y}_1^n|\overline{\mathbf{Y}}_0) = \int p_\theta(\mathbf{Y}_1^n|\overline{\mathbf{Y}}_0, X_0 = x_0)\xi(dx_0).$$
That is, instead of choosing an initial point $X_0 = x_0$ we set instead an arbitrary initial distribution. There is in general little rationale for doing that, but the results obtained below for a fixed initial condition $X_0 = x_0$ immediately carry over to this more general context. Typically such a $\xi$ has a density w.r.t. $\mu$ so that there are a density $p_{\theta,\xi}(\mathbf{Y}_1^n|\overline{\mathbf{Y}}_0, X_0 = x_0)$ and an associated MLE.

The consistency proof for the MLE follows the classical scheme of Wald (1949), which amounts to proving that there exists a deterministic asymptotic criterion function $l(\theta)$ such that $n^{-1}l_n(\theta, x_0) \to l(\theta)$ $\overline{\mathbb{P}}_{\theta^*}$-a.s. uniformly w.r.t. $\theta \in \Theta$ and that $\theta^*$ is a well-separated point of maximum of $l(\theta)$. It should be stressed that the asymptotic criterion $l(\theta)$ should of course not depend on the initial point $X_0 = x_0$.

The first step of the proof thus consists in showing that the normalized log likelihood function $n^{-1}l_n(\theta, x_0)$ converges to $l(\theta)$ uniformly w.r.t. $\theta$. This requires a uniform (w.r.t. $\theta \in \Theta$ and $x_0 \in \mathcal{X}$) law of large numbers. We first show that the difference between the conditional log likelihood functions $l_n(\theta, x_0)$ and $l_n(\theta)$ stays within some deterministic bound, and hence $n^{-1}(l_n(\theta, x_0) - l_n(\theta))$ tends to zero $\overline{\mathbb{P}}_{\theta^*}$-a.s. and in $L^1(\overline{\mathbb{P}}_{\theta^*})$ [see Del Moral and Miclo (2001) for similar results].

LEMMA 2. *Assume* (A1) *and* (A2). *Then, for all* $x_0 \in \mathcal{X}$,
$$\sup_{\theta \in \Theta} |l_n(\theta, x_0) - l_n(\theta)| \le 1/(1-\rho)^2, \qquad \overline{\mathbb{P}}_{\theta^*}\text{-}a.s.$$

PROOF. Note that by Corollary 1, (8) and (11),
$$\|\mathbb{P}_\theta(X_{k-1} \in \cdot|\overline{\mathbf{Y}}_0^{k-1}, X_0 = x_0) - \overline{\mathbb{P}}_\theta(X_{k-1} \in \cdot|\overline{\mathbf{Y}}_0^{k-1})\|_{\mathrm{TV}} \le \rho^{k-1}.$$
This implies that, for $k \ge 1$,
$$|p_\theta(Y_k|\overline{\mathbf{Y}}_0^{k-1}, X_0 = x_0) - \bar{p}_\theta(Y_k|\overline{\mathbf{Y}}_0^{k-1})|$$
$$= \left|\int\int g_\theta(Y_k|\overline{\mathbf{Y}}_{k-1}, x_k) q_\theta(x_{k-1}, x_k) \mu(dx_k) \right.$$
$$\left. \times (\mathbb{P}_\theta(dx_{k-1}|\overline{\mathbf{Y}}_0^{k-1}, X_0 = x_0) - \overline{\mathbb{P}}_\theta(dx_{k-1}|\overline{\mathbf{Y}}_0^{k-1}))\right|$$
$$\le \rho^{k-1} \sup_{x_{k-1}} \int g_\theta(Y_k|\overline{\mathbf{Y}}_{k-1}, x_k) q_\theta(x_{k-1}, x_k) \mu(dx_k)$$
$$\le \rho^{k-1} \sigma_+ \int g_\theta(Y_k|\overline{\mathbf{Y}}_{k-1}, x) \mu(dx).$$

In addition, by (7),

$$p_\theta(Y_k|\overline{\mathbf{Y}}_0^{k-1}, X_0 = x_0) \geq \sigma_- \int g_\theta(Y_k|\overline{\mathbf{Y}}_{k-1}, x)\mu(dx),$$

and the same inequality holds for $\bar{p}_\theta(Y_k|\overline{\mathbf{Y}}_0^{k-1})$. The inequality $|\log x - \log y| \leq |x - y|/(x \wedge y)$ now shows that

$$|\log p_\theta(Y_k|\overline{\mathbf{Y}}_0^{k-1}, X_0 = x_0) - \log \bar{p}_\theta(Y_k|\overline{\mathbf{Y}}_0^{k-1})| \leq \frac{\rho^{k-1}}{1-\rho}.$$

A summation concludes the proof. □

The next step consists in showing that $n^{-1}l_n(\theta)$ can be approximated by the sample mean of a $\overline{\mathbb{P}}_{\theta^*}$-stationary ergodic sequence of random variables in $L^1(\mathbb{P}_{\theta^*})$. It is natural to approximate $n^{-1}l_n(\theta) = n^{-1}\sum_{k=1}^n \log \bar{p}_\theta(Y_k|\overline{\mathbf{Y}}_0^{k-1})$ by $n^{-1}\sum_{k=1}^n \log \bar{p}_\theta(Y_k|\overline{\mathbf{Y}}_{-\infty}^{k-1})$, provided we can give meaning to the latter conditional densities. This is the main purpose of the construction below. Let, for $x \in \mathcal{X}$,

$$\Delta_{k,m,x}(\theta) \triangleq \log \bar{p}_\theta(Y_k|\overline{\mathbf{Y}}_{-m}^{k-1}, X_{-m} = x),$$

$$\Delta_{k,m}(\theta) \triangleq \log \bar{p}_\theta(Y_k|\overline{\mathbf{Y}}_{-m}^{k-1}) = \log \int \bar{p}_\theta(Y_k|\overline{\mathbf{Y}}_{-m}^{k-1}, X_{-m} = x_{-m})\overline{\mathbb{P}}_\theta(dx_{-m}|\overline{\mathbf{Y}}_{-m}^{k-1}).$$

It follows from the definitions that $l_n(\theta) = \sum_{k=1}^n \Delta_{k,0}(\theta)$. In order to show that for any $k \geq 0$ the sequences $\{\Delta_{k,m}(\theta)\}_{m\geq 0}$ and $\{\Delta_{k,m,x}(\theta)\}_{m\geq 0}$ converge uniformly w.r.t. $\theta \in \Theta$, $\overline{\mathbb{P}}_{\theta^*}$-a.s., we prove that they are uniform Cauchy sequences. This property is implied by the following lemma.

LEMMA 3. *Assume* (A1)–(A3). *Then for all $k \geq 1$ and $m, m' \geq 0$, $\overline{\mathbb{P}}_{\theta^*}$-a.s.,*

(12) $\quad \sup_{\theta \in \Theta} \sup_{x,x' \in \mathcal{X}} |\Delta_{k,m,x}(\theta) - \Delta_{k,m',x'}(\theta)| \leq \rho^{k+(m\wedge m')-1}/(1-\rho),$

(13) $\quad \sup_{\theta \in \Theta} \sup_{x \in \mathcal{X}} |\Delta_{k,m,x}(\theta) - \Delta_{k,m}(\theta)| \leq \rho^{k+m-1}/(1-\rho),$

(14) $\quad \sup_{\theta \in \Theta} \sup_{m \geq 0} \sup_{x \in \mathcal{X}} |\Delta_{k,m,x}(\theta)| \leq \max\left(|\log b_+|, |\log (\sigma_- b_-(\overline{\mathbf{Y}}_{k-1}, Y_k))|\right).$

The proof is similar to the proof of Lemma 2 (making use of the uniform ergodicity of the conditional chain) and is given in the Appendix. By (12), $\{\Delta_{k,m,x}(\theta)\}_{m\geq 0}$ is a uniform Cauchy sequence w.r.t. $\theta \in \Theta$ $\overline{\mathbb{P}}_{\theta^*}$-a.s. and thus $\Delta_{k,m,x}(\theta)$ converges uniformly $\overline{\mathbb{P}}_{\theta^*}$-a.s. Equation (12) also implies that $\lim_{m\to\infty} \Delta_{k,m,x}(\theta)$ does not depend on $x$. Denote by $\Delta_{k,\infty}(\theta)$ this limit. Intuitively, we may think of $\Delta_{k,\infty}(\theta)$ as $\log \bar{p}_\theta(Y_k|\mathbf{Y}_{-\infty}^{k-1})$. Equation (14) shows that $\{\Delta_{k,m,x}(\theta)\}_{m\geq 0}$ is uniformly bounded in $L^1(\overline{\mathbb{P}}_{\theta^*})$, and thus the limit $\Delta_{k,\infty}(\theta)$ is also in $L^1(\overline{\mathbb{P}}_{\theta^*})$. Note that $\{\Delta_{k,\infty}(\theta)\}$ is a $\overline{\mathbb{P}}_{\theta^*}$-stationary ergodic process.

Setting $m = 0$ in (12) and letting $m' \to \infty$ shows that, $\overline{\mathbb{P}}_{\theta^*}$-a.s.,

$$\sup_{\theta \in \Theta} |\Delta_{k,0,x}(\theta) - \Delta_{k,\infty}(\theta)| \leq \rho^{k-1}/(1-\rho),$$

and (13) shows that $\sup_{\theta \in \Theta} |\Delta_{k,0,x}(\theta) - \Delta_{k,0}(\theta)| \leq \rho^{k-1}/(1-\rho)$. These two relations readily imply the following result.

COROLLARY 2. *Assume* (A1) *and* (A2). *Then*

$$\sum_{k=1}^n \sup_{\theta \in \Theta} |\Delta_{k,0}(\theta) - \Delta_{k,\infty}(\theta)| \leq 2/(1-\rho)^2, \qquad \overline{\mathbb{P}}_{\theta^*}\text{-a.s.}$$

Corollary 2 shows that $n^{-1}l_n(\theta)$ can be approximated by the sample mean of a stationary ergodic sequence, uniformly w.r.t. $\theta \in \Theta$. Since $\Delta_{0,\infty}(\theta) \in L^1(\overline{\mathbb{P}}_{\theta^*})$, the ergodic theorem implies that $n^{-1}l_n(\theta) \to l(\theta) \triangleq \overline{\mathbb{E}}_{\theta^*}[\Delta_{0,\infty}(\theta)]$ $\overline{\mathbb{P}}_{\theta^*}$-a.s. and in $L^1(\overline{\mathbb{P}}_{\theta^*})$. Combining this result with Lemma 2 yields the following.

PROPOSITION 1. *Assume* (A1)–(A3). *Then for all* $x_0 \in \mathcal{X}$ *and* $\theta \in \Theta$,

$$\lim_{n \to \infty} n^{-1}l_n(\theta, x_0) = l(\theta), \qquad \overline{\mathbb{P}}_{\theta^*}\text{-a.s. and in } L^1(\overline{\mathbb{P}}_{\theta^*}).$$

REMARK 4. The pointwise convergence of $n^{-1}l_n(\theta, x_0)$ has been established for HMMs by Leroux (1992) and Le Gland and Mevel (2000) for a finite state space and later for a compact state space by Douc and Matias (2001). In the papers of Le Gland and Mevel (2000) and Douc and Matias (2001), the authors used the geometric ergodicity of an extended Markov chain consisting of the hidden variable, the observed variable and the prediction filter density function. However, this approach requires conditions stronger than the weak ergodicity condition (A2) and the moment condition (A3).

The next step of the proof consists in showing that $l(\theta)$ is continuous w.r.t. $\theta$. To that purpose, first observe that, by (14) and the dominated convergence theorem, for any $x \in \mathcal{X}$ and $\theta \in \Theta$,

$$l(\theta) = \overline{\mathbb{E}}_{\theta^*}\left[\lim_{m \to \infty} \Delta_{0,m,x}(\theta)\right] = \lim_{m \to \infty} \overline{\mathbb{E}}_{\theta^*}[\Delta_{0,m,x}(\theta)].$$

Since $\{\Delta_{0,m,x}(\theta)\}_{m\geq 0}$ is a uniform Cauchy sequence $\overline{\mathbb{P}}_{\theta^*}$-a.s. which is uniformly bounded in $L^1(\overline{\mathbb{P}}_{\theta^*})$ ($\overline{\mathbb{E}}_{\theta^*}[\sup_{m\geq 0} \sup_{\theta \in \Theta}[\Delta_{0,m,x}(\theta)]] < \infty$), it suffices to show that $\Delta_{0,m,x}(\theta)$ is continuous w.r.t. $\theta$. In fact, this is the whole point of using $\Delta_{0,m,x}(\theta)$ instead of $\Delta_{0,m}(\theta)$. We will need the following additional assumption:

(A4) For all $x, x' \in \mathcal{X}$ and all $(\bar{\mathbf{y}}, y') \in \mathcal{Y}^s \times \mathcal{Y}$, $\theta \mapsto q_\theta(x, x')$ and $\theta \mapsto g_\theta(y'|\bar{\mathbf{y}}, x)$ are continuous.

LEMMA 4. *Assume* (A1)–(A4). *Then for all* $\theta \in \Theta$,
$$\lim_{\delta \to 0} \overline{\mathbb{E}}_{\theta^*} \left[ \sup_{|\theta' - \theta| \leq \delta} |\Delta_{0,\infty}(\theta') - \Delta_{0,\infty}(\theta)| \right] = 0.$$

The proof is given in the Appendix. We may now state the central result of this section, the uniform convergence of the normalized log likelihood $n^{-1} l_n(\theta, x_0)$ to $l(\theta)$, which follows almost immediately from Corollary 2 and the ergodic theorem.

PROPOSITION 2. *Assume* (A1)–(A4). *Then*
$$\lim_{n \to \infty} \sup_{\theta \in \Theta} \sup_{x_0 \in \mathcal{X}} |n^{-1} l_n(\theta, x_0) - l(\theta)| = 0, \qquad \overline{\mathbb{P}}_{\theta^*}\text{-}a.s.$$

Again, the proof is in the Appendix.

**5. Consistency of the maximum likelihood estimator.** We will now prove that under suitable assumptions the unique maximizer of $\theta \mapsto l(\theta)$ is $\theta^*$, the true value of the parameter. Let $\overline{\mathbb{P}}_\theta^Y$ be the trace of $\overline{\mathbb{P}}_\theta$ on $\{\mathcal{Y}^\mathbb{N}, \mathcal{B}(\mathcal{Y})^{\otimes \mathbb{N}}\}$, that is, the distribution of $\{Y_k\}$. Consider the following assumption:

(A5) $\theta = \theta^*$ if and only if

(15) $$\overline{\mathbb{P}}_\theta^Y = \overline{\mathbb{P}}_{\theta^*}^Y.$$

In other words, under this assumption the stationary laws of the observed process associated with two different values of the parameter do not coincide unless the parameters do. This is obviously the minimal assumption that we can impose. When it comes to applying the results, it is sometimes more convenient to consider the following alternative identifiability condition, which in certain circumstances proves easier to verify.

(A5′) $\theta = \theta^*$ if and only if

(16) $$\overline{\mathbb{E}}_{\theta^*} \left[ \log \frac{\bar{p}_{\theta^*}(\mathbf{Y}_1^p | \overline{\mathbf{Y}}_0)}{\bar{p}_\theta(\mathbf{Y}_1^p | \overline{\mathbf{Y}}_0)} \right] = 0 \qquad \text{for all } p \geq 1.$$

In fact we will see below that under (A1)–(A3), these two conditions are equivalent. Of course neither of the identifiability assumptions stated above is entirely satisfactory, because both conditions implicitly make use of the stationary distribution of the complete chain, which typically is infeasible to compute. Nevertheless, it does not seem sensible to expect much simpler identifiability conditions based, say, on $g_\theta$ and $q_\theta$ alone. The usefulness of (A5′) is revealed when conditioning on $\overline{\mathbf{Y}}_0$, yielding $\theta = \theta^*$ if and only if

$$(17) \qquad \overline{\mathbb{E}}_{\theta^*}\left[\overline{\mathbb{E}}_{\theta^*}\left(\log \frac{\bar{p}_{\theta^*}(\mathbf{Y}_1^p|\overline{\mathbf{Y}}_0)}{\bar{p}_{\theta}(\mathbf{Y}_1^p|\overline{\mathbf{Y}}_0)}\Big|\overline{\mathbf{Y}}_0\right)\right] = 0 \qquad \text{for all } p \geq 1.$$

In this expression the inner expectation is a conditional Kullback–Leibler measure, and hence nonnegative. If equality holds in (17), the inner conditional expectation vanishes $\overline{\mathbb{P}}_{\theta^*}^{\overline{\mathbf{Y}}_0}$-a.s. This observation may in turn often be used to prove that $\theta = \theta^*$, using, for example, identifiability of mixtures of the family to which the densities $g_\theta(\cdot|\bar{\mathbf{y}}, x)$ belong. A particular example involving linear regressions with normal disturbances and finite-valued regime is discussed in Krishnamurthy and Rydén [(1998), page 302]. Slightly different identifiability conditions are employed in Francq and Roussignol (1998).

Before proceeding to the equivalence of (A5) and (A5′), some preparatory lemmas are needed. We will first show that the conditional density function $\bar{p}_\theta(\mathbf{Y}_k^\ell|\overline{\mathbf{Y}}_i^j)$ $(i \leq j < k \leq \ell)$ converges to the unconditional density function $\bar{p}_\theta(\mathbf{Y}_k^\ell)$ when the gap $k-j$ tends to infinity. This can be viewed as a kind of $\phi$-mixing condition expressed directly on the conditional and the unconditional density functions, which is inherited from the ergodicity of the complete chain.

LEMMA 5. *Assume* (A1)–(A3) *and fix* $k \leq \ell$. *Then*

$$\lim_{j \to -\infty} \sup_{i \leq j} |\bar{p}_\theta(\mathbf{Y}_k^\ell|\overline{\mathbf{Y}}_i^j) - \bar{p}_\theta(\mathbf{Y}_k^\ell)| = 0 \qquad \text{in } \overline{\mathbb{P}}_{\theta^*}\text{-probability.}$$

The proof is given in the Appendix.

The following lemma shows that (15) and (16) are equivalent.

LEMMA 6. *Assume* (A1)–(A3). *Then* (15) *holds if and only if* (16) *holds.*

PROOF. It obviously suffices to show the "if" part, so suppose (16) holds. The basic idea consists in inserting a gap in the range of variables. For $p \geq 1$ and $m \geq 0$, write

$$0 = \overline{\mathbb{E}}_{\theta^*}\left[\log \frac{\bar{p}_{\theta^*}(\mathbf{Y}_1^{p+m}|\overline{\mathbf{Y}}_0)}{\bar{p}_{\theta}(\mathbf{Y}_1^{p+m}|\overline{\mathbf{Y}}_0)}\right]$$

$$= \overline{\mathbb{E}}_{\theta^*}\left[\log \frac{\bar{p}_{\theta^*}(\mathbf{Y}_1^m|\mathbf{Y}_{m+1}^{p+m}, \overline{\mathbf{Y}}_0)}{\bar{p}_{\theta}(\mathbf{Y}_1^m|\mathbf{Y}_{m+1}^{p+m}, \overline{\mathbf{Y}}_0)}\right] + \overline{\mathbb{E}}_{\theta^*}\left[\log \frac{\bar{p}_{\theta^*}(\mathbf{Y}_{m+1}^{p+m}|\overline{\mathbf{Y}}_0)}{\bar{p}_{\theta}(\mathbf{Y}_{m+1}^{p+m}|\overline{\mathbf{Y}}_0)}\right].$$

The two terms on the right-hand side are expectations of Kullback–Leibler divergence functions and thus nonnegative, which shows that

$$0 \geq \overline{\mathbb{E}}_{\theta^*}\left[\log \frac{\bar{p}_{\theta^*}(\mathbf{Y}_{m+1}^{p+m}|\overline{\mathbf{Y}}_0)}{\bar{p}_{\theta}(\mathbf{Y}_{m+1}^{p+m}|\overline{\mathbf{Y}}_0)}\right] = \overline{\mathbb{E}}_{\theta^*}\left[\log \frac{\bar{p}_{\theta^*}(\mathbf{Y}_1^p|\overline{\mathbf{Y}}_{-m})}{\bar{p}_{\theta}(\mathbf{Y}_1^p|\overline{\mathbf{Y}}_{-m})}\right]$$

$$= \overline{\mathbb{E}}_{\theta^*}\left[\int \log \frac{\bar{p}_{\theta^*}(\mathbf{Y}_1^p = \mathbf{y}_1^p|\overline{\mathbf{Y}}_{-m})}{\bar{p}_{\theta}(\mathbf{Y}_1^p = \mathbf{y}_1^p|\overline{\mathbf{Y}}_{-m})}\bar{p}_{\theta^*}(\mathbf{Y}_1^p = \mathbf{y}_1^p|\overline{\mathbf{Y}}_{-m})\nu^{\otimes p}(d\mathbf{y}_1^p)\right].$$

Thus, for all $m \geq 0$,

$$\bar{p}_{\theta^*}(\mathbf{Y}_1^p|\overline{\mathbf{Y}}_{-m}) = \bar{p}_{\theta}(\mathbf{Y}_1^p|\overline{\mathbf{Y}}_{-m}), \qquad \overline{\mathbb{P}}_{\theta^*}\text{-a.s.}$$

By Lemma 5,

$$|\bar{p}_{\theta^*}(\mathbf{Y}_1^p) - \bar{p}_{\theta}(\mathbf{Y}_1^p)|$$
$$= \lim_{m \to \infty} |\bar{p}_{\theta^*}(\mathbf{Y}_1^p|\overline{\mathbf{Y}}_{-m}) - \bar{p}_{\theta}(\mathbf{Y}_1^p|\overline{\mathbf{Y}}_{-m})| = 0 \qquad \text{in } \overline{\mathbb{P}}_{\theta^*}\text{-probability},$$

whence $\bar{p}_{\theta^*}(\mathbf{Y}_1^p) = \bar{p}_{\theta}(\mathbf{Y}_1^p)$ $\overline{\mathbb{P}}_{\theta^*}$-a.s. The proof is complete. $\square$

PROPOSITION 3. *Under* (A1)–(A5), $l(\theta) \leq l(\theta^*)$ *and* $l(\theta) = l(\theta^*)$ *if and only if* $\theta = \theta^*$.

PROOF. By the dominated convergence theorem,

$$l(\theta) = \overline{\mathbb{E}}_{\theta^*}\left[\lim_{m \to \infty} \log \bar{p}_{\theta}(Y_1|\overline{\mathbf{Y}}_{-m}^0)\right]$$
(18)
$$= \lim_{m \to \infty} \overline{\mathbb{E}}_{\theta^*}[\log \bar{p}_{\theta}(Y_1|\overline{\mathbf{Y}}_{-m}^0)]$$
$$= \lim_{m \to \infty} \overline{\mathbb{E}}_{\theta^*}[\overline{\mathbb{E}}_{\theta^*}[\log \bar{p}_{\theta}(Y_1|\overline{\mathbf{Y}}_{-m}^0)|\overline{\mathbf{Y}}_{-m}^0]].$$

Hence $l(\theta^*) - l(\theta)$ is nonnegative as the limit of expectations of conditional Kullback–Leibler divergence functions and $\theta^*$ is a maximizer of the function $\theta \mapsto l(\theta)$.

Now assume $l(\theta) = l(\theta^*)$. By Lemma 6 it suffices to prove that (16) holds. Note that for any $k \geq 1$ and $m \geq 0$,

$$\overline{\mathbb{E}}_{\theta^*}[\log \bar{p}_{\theta}(\mathbf{Y}_1^k|\overline{\mathbf{Y}}_{-m}^0)] = \sum_{i=1}^{k} \overline{\mathbb{E}}_{\theta^*}[\log \bar{p}_{\theta}(Y_1|\overline{\mathbf{Y}}_{-m-i+1}^0)].$$

Hence, by (18),

$$\lim_{m \to \infty} \overline{\mathbb{E}}_{\theta^*}[\log \bar{p}_{\theta}(\mathbf{Y}_1^k|\overline{\mathbf{Y}}_{-m}^0)] = kl(\theta),$$

and for $p + s < k + 1$,

$$\begin{aligned}
0 = k(l(\theta^*) - l(\theta)) &= \lim_{m \to \infty} \overline{\mathbb{E}}_{\theta^*} \left[ \log \frac{\bar{p}_{\theta^*}(\mathbf{Y}_1^k | \overline{\mathbf{Y}}_{-m}^0)}{\bar{p}_{\theta}(\mathbf{Y}_1^k | \overline{\mathbf{Y}}_{-m}^0)} \right] \\
&\geq \limsup_{m \to \infty} \overline{\mathbb{E}}_{\theta^*} \left[ \log \frac{\bar{p}_{\theta^*}(\mathbf{Y}_{k-p+1}^k | \overline{\mathbf{Y}}_{k-p}, \overline{\mathbf{Y}}_{-m}^0)}{\bar{p}_{\theta}(\mathbf{Y}_{k-p+1}^k | \overline{\mathbf{Y}}_{k-p}, \overline{\mathbf{Y}}_{-m}^0)} \right] \\
&= \limsup_{m \to \infty} \overline{\mathbb{E}}_{\theta^*} \left[ \log \frac{\bar{p}_{\theta^*}(\mathbf{Y}_1^p | \overline{\mathbf{Y}}_0, \overline{\mathbf{Y}}_{p-k-m}^{p-k})}{\bar{p}_{\theta}(\mathbf{Y}_1^p | \overline{\mathbf{Y}}_0, \overline{\mathbf{Y}}_{p-k-m}^{p-k})} \right].
\end{aligned}$$

The proof is concluded by letting $k \to \infty$ and using Lemma 7. □

LEMMA 7. *Assume* (A1)–(A3). *Then for all* $p \geq 1$ *and all* $\theta \in \Theta$,

$$\lim_{k \to \infty} \sup_{m \geq k} \left| \overline{\mathbb{E}}_{\theta^*} \left[ \log \frac{\bar{p}_{\theta^*}(\mathbf{Y}_1^p | \overline{\mathbf{Y}}_0, \overline{\mathbf{Y}}_{-m}^{-k})}{\bar{p}_{\theta}(\mathbf{Y}_1^p | \overline{\mathbf{Y}}_0, \overline{\mathbf{Y}}_{-m}^{-k})} \right] - \overline{\mathbb{E}}_{\theta^*} \left[ \log \frac{\bar{p}_{\theta^*}(\mathbf{Y}_1^p | \overline{\mathbf{Y}}_0)}{\bar{p}_{\theta}(\mathbf{Y}_1^p | \overline{\mathbf{Y}}_0)} \right] \right| = 0.$$

The proof of this lemma is based on the mixing properties of the complete chain (see Lemma 5) and is postponed to the Appendix.

We may now summarize our findings in the following theorem, which states the strong consistency of the conditional MLE.

THEOREM 1. *Assume* (A1)–(A5). *Then, for any* $x_0 \in \mathcal{X}$, $\lim_{n \to \infty} \hat{\theta}_{n,x_0} = \theta^*$, $\overline{\mathbb{P}}_{\theta^*}$-*a.s.*

**6. Asymptotic normality of the maximum likelihood estimator.** Lemma 1 and Corollary 1 are the basic tools for generalizing the results of BRR and JP. The pattern of the proof of asymptotic normality of the MLE is similar to that presented in these contributions, with two major differences. First, the geometric upper bounds are deterministic, which is a consequence of Lemma 1 and Corollary 1. Second, in this paper, the MLE is the maximizer of the conditional log likelihood $l_n(\theta, x_0)$, where $x_0$ is some fixed arbitrary point in $\mathcal{X}$, whereas in BRR and JP it is the maximizer of the unconditional log likelihood $l_n(\theta)$.

Not surprisingly, the proof of asymptotic normality requires some additional regularity assumptions. Let $\nabla_\theta$ and $\nabla_\theta^2$ be the gradient and the Hessian operator with respect to the parameter $\theta$, respectively. We will assume that there exists a positive real $\delta$ such that on $G \triangleq \{\theta \in \Theta : |\theta - \theta^*| < \delta\}$ the following conditions hold:

(A6) For all $x, x' \in \mathcal{X}$ and $(\bar{\mathbf{y}}, y') \in \mathcal{Y}^s \times \mathcal{Y}$, the functions $\theta \mapsto q_\theta(x, x')$ and $\theta \mapsto g_\theta(y' | \bar{\mathbf{y}}, x')$ are twice continuously differentiable on $G$.

(A7) (a) $\sup_{\theta\in G}\sup_{x,x'}\|\nabla_\theta\log q_\theta(x,x')\|<\infty$ and $\sup_{\theta\in G}\sup_{x,x'}\|\nabla_\theta^2\times\log q_\theta(x,x')\|<\infty$.

(b) $\overline{\mathbb{E}}_{\theta^*}[\sup_{\theta\in G}\sup_x\|\nabla_\theta\log g_\theta(Y_1|\overline{\mathbf{Y}}_0,x)\|^2]<\infty$ and $\overline{\mathbb{E}}_{\theta^*}[\sup_{\theta\in G}\sup_x\|\nabla_\theta^2\log g_\theta(Y_1|\overline{\mathbf{Y}}_0,x)\|]<\infty$.

(A8) (a) For $\bar{\nu}\otimes\nu$-almost all $(\bar{\mathbf{y}},y')$ in $\mathcal{Y}^s\times\mathcal{Y}$ there exists a function $f_{\bar{\mathbf{y}},y'}:\mathcal{X}\to\mathbb{R}^+$ in $L^1(\mu)$ such that $\sup_{\theta\in G}g_\theta(y'|\bar{\mathbf{y}},x)\leq f_{\bar{\mathbf{y}},y'}(x)$.

(b) For $\mu\otimes\bar{\nu}$-almost all $(x,\bar{\mathbf{y}})\in\mathcal{X}\times\mathcal{Y}^s$, there exist functions $f^1_{x,\bar{\mathbf{y}}}:\mathcal{Y}\to\mathbb{R}^+$ and $f^2_{x,\bar{\mathbf{y}}}:\mathcal{Y}\to\mathbb{R}^+$ in $L^1(\nu)$ such that $\|\nabla_\theta g_\theta(y'|\bar{\mathbf{y}},x)\|\leq f^1_{x,\bar{\mathbf{y}}}(y')$ and $\|\nabla_\theta^2 g_\theta(y'|\bar{\mathbf{y}},x)\|\leq f^2_{x,\bar{\mathbf{y}}}(y')$ for all $\theta\in G$.

REMARK 5. The regularity requirements (existence of first and second derivatives at all points, existence of integrable upper bounds) are reminiscent of Cramér's classical proof of asymptotic normality of the MLE. It is obvious that these conditions could have been weakened using more sophisticated techniques. We will nevertheless stick to the conventional proof.

REMARK 6. The conditions are weaker and more easily checked than those used by JP, who assumed that the stationary density of the complete Markov chain is twice differentiable w.r.t. to $\theta$, a condition which is difficult to check except for very simple models. However, as seen below, by using proper conditioning techniques it is possible to avoid such assumptions.

Asymptotic normality of the MLE is implied by:

(i) a central limit theorem (CLT) for the Fisher score function $n^{-1/2}\nabla_\theta l_n(\theta^*,x_0)$, and

(ii) a locally uniform law of large numbers for the observed Fisher information $-n^{-1}\nabla_\theta^2 l_n(\theta,x_0)$ for $\theta$ in a neighborhood of $\theta^*$.

Along the lines of the proofs by BRR and JP, the key to the proof consists in finding proper expressions for these two quantities. Exploiting the hierarchical structure of the model, it turns out that it is practical to express the score function and the observed Fisher information as functions of conditional expectations of the complete score function and the complete Fisher information.

6.1. *A central limit theorem for the score function.* The Fisher identity [Louis (1982)] generally states that for a model with missing data, the score function equals the conditional expectation of the complete score given the observed data; the complete score is the gradient of the complete log likelihood, that is, the likelihood that includes the missing data in addition to the observed data. The rationale for using this identity is that the log likelihood and score functions themselves are typically rather involved [cf.

(1)] while the complete log likelihood and score are simpler. This is true in our case, in which the Markov chain $\{X_k\}_{k=1}^n$ constitutes the missing data. The Fisher identity requires exchanging the gradient operator with certain integrals, and is valid under (A7) and (A8). Hence, for any $x_0 \in \mathcal{X}$,

$$n^{-1/2}\nabla_\theta l_n(\theta^*, x_0) = n^{-1/2} \sum_{k=1}^n \nabla_\theta \log p_{\theta^*}(Y_k|\overline{\mathbf{Y}}_0^{k-1}, X_0 = x_0)$$

$$= n^{-1/2} \sum_{k=1}^n \Delta_{k,0,x_0}(\theta^*),$$

where for any $x \in \mathcal{X}$ and $\theta \in \Theta$,

$$\Delta_{k,0,x}(\theta) \triangleq \mathbb{E}_\theta\left[\sum_{i=1}^k \phi(\theta, \overline{\mathbf{Z}}_{i-1}^i)\Big|\overline{\mathbf{Y}}_0^k, X_0 = x\right]$$
$$- \mathbb{E}_\theta\left[\sum_{i=1}^{k-1} \phi(\theta, \overline{\mathbf{Z}}_{i-1}^i)\Big|\overline{\mathbf{Y}}_0^{k-1}, X_0 = x\right],$$

with the convention $\sum_{i=a}^b c_i = 0$ if $b < a$ and

$$\phi(\theta, \overline{\mathbf{Z}}_{i-1}^i) = \phi(\theta, \overline{\mathbf{Z}}_{i-1}, \overline{\mathbf{Z}}_i) = \phi(\theta, (\mathbf{X}_{i-s}^{i-1}, \overline{\mathbf{Y}}_{i-1}), (\mathbf{X}_{i-s+1}^i, \mathbf{Y}_{i-s+1}^i))$$
$$= \phi(\theta, (X_{i-1}, \overline{\mathbf{Y}}_{i-1}), (X_i, Y_i))$$
$$\triangleq \nabla_\theta \log\left(q_\theta(X_{i-1}, X_i) g_\theta(Y_i|\overline{\mathbf{Y}}_{i-1}, X_i)\right)$$

is the conditional score of $(X_i, Y_i)$ given $(X_{i-1}, \overline{\mathbf{Y}}_{i-1})$.

We also let, for $m \geq 0$,

$$\Delta_{k,m}(\theta) \triangleq \overline{\mathbb{E}}_\theta\left[\sum_{i=-m+1}^k \phi(\theta, \overline{\mathbf{Z}}_{i-1}^i)\Big|\overline{\mathbf{Y}}_{-m}^k\right] - \overline{\mathbb{E}}_\theta\left[\sum_{i=-m+1}^{k-1} \phi(\theta, \overline{\mathbf{Z}}_{i-1}^i)\Big|\overline{\mathbf{Y}}_{-m}^{k-1}\right].$$

Similar to what is done in BRR and JP we show that $\Delta_{k,0,x}(\theta^*)$ can be approximated in $L^2(\overline{\mathbb{P}}_{\theta^*})$ by a $\overline{\mathbb{P}}_{\theta^*}$-stationary martingale increment sequence and apply a CLT for sums of stationary martingale increments.

The first step in the proof consists of showing that the initial point $x$ does not show up in the limit.

LEMMA 8. *Assume* (A1), (A2) *and* (A6)–(A7). *Then, for all* $x \in \mathcal{X}$,

$$\lim_{n \to \infty} \overline{\mathbb{E}}_{\theta^*}\left\|n^{-1/2} \sum_{k=1}^n (\Delta_{k,0,x}(\theta^*) - \Delta_{k,0}(\theta^*))\right\|^2 = 0.$$

PROOF. Write
$$\sum_{k=1}^{n} (\Delta_{k,0,x}(\theta^*) - \Delta_{k,0}(\theta^*))$$
$$= \sum_{k=1}^{n} (\mathbb{E}_{\theta^*}[\phi(\theta^*, \overline{\mathbf{Z}}_{k-1}^k)|\overline{\mathbf{Y}}_0^n, X_0 = x] - \overline{\mathbb{E}}_{\theta^*}[\phi(\theta^*, \overline{\mathbf{Z}}_{k-1}^k)|\overline{\mathbf{Y}}_0^n]).$$

Under the stated assumptions $\overline{\mathbb{E}}_{\theta^*}(\sup_{x,x' \in \mathcal{X}} \|\phi(\theta^*, (x, \overline{\mathbf{Y}}_0), (x', Y_1))\|^2) < \infty$. The proof now follows from Corollary 1, which implies that
$$\|\mathbb{E}_{\theta^*}[\phi(\theta^*, \overline{\mathbf{Z}}_{k-1}^k)|\overline{\mathbf{Y}}_0^n, X_0 = x] - \overline{\mathbb{E}}_{\theta^*}[\phi(\theta^*, \overline{\mathbf{Z}}_{k-1}^k)|\overline{\mathbf{Y}}_0^n]\|$$
$$\leq 2 \sup_{x,x' \in \mathcal{X}} \|\phi(\theta^*, (x, \overline{\mathbf{Y}}_{k-1}), (x', Y_k))\| \rho^{k-1}. \qquad \Box$$

We will now show that for any $k$, $\{\Delta_{k,m}(\theta^*)\}_{m \geq 0}$ is a Cauchy sequence in $L^2(\overline{\mathbb{P}}_{\theta^*})$. Since
$$\Delta_{k,m}(\theta^*) = \overline{\mathbb{E}}_{\theta^*}[\phi(\theta^*, \overline{\mathbf{Z}}_{k-1}^k)|\overline{\mathbf{Y}}_{-m}^k]$$
$$+ \sum_{i=-m+1}^{k-1} (\overline{\mathbb{E}}_{\theta^*}[\phi(\theta^*, \overline{\mathbf{Z}}_{i-1}^i)|\overline{\mathbf{Y}}_{-m}^k] - \overline{\mathbb{E}}_{\theta^*}[\phi(\theta^*, \overline{\mathbf{Z}}_{i-1}^i)|\overline{\mathbf{Y}}_{-m}^{k-1}]),$$

the difference $\Delta_{k,m}(\theta) - \Delta_{k,m'}(\theta)$ (assuming $m' > m > 0$) involves for each $-m < i \leq k$ terms of the form either $\overline{\mathbb{E}}_{\theta^*}[\phi(\theta^*, \overline{\mathbf{Z}}_{i-1}^i)|\overline{\mathbf{Y}}_{-m}^k] - \overline{\mathbb{E}}_{\theta^*}[\phi(\theta^*, \overline{\mathbf{Z}}_{i-1}^i)|\overline{\mathbf{Y}}_{-m'}^k]$ or $\overline{\mathbb{E}}_{\theta^*}[\phi(\theta^*, \overline{\mathbf{Z}}_{i-1}^i)|\overline{\mathbf{Y}}_{-m}^k] - \overline{\mathbb{E}}_{\theta^*}[\phi(\theta^*, \overline{\mathbf{Z}}_{i-1}^i)|\overline{\mathbf{Y}}_{-m}^{k-1}]$. By Corollary 1 and an argument used to prove Lemma 3 we obtain for $-m' < -m < i \leq k$ that, $\overline{\mathbb{P}}_{\theta^*}$-a.s.,

$$\|\overline{\mathbb{E}}_{\theta^*}[\phi(\theta^*, \overline{\mathbf{Z}}_{i-1}^i)|\overline{\mathbf{Y}}_{-m}^k] - \overline{\mathbb{E}}_{\theta^*}[\phi(\theta^*, \overline{\mathbf{Z}}_{i-1}^i)|\overline{\mathbf{Y}}_{-m'}^k]\|$$
(19)
$$\leq 2 \sup_{x,x' \in \mathcal{X}} \|\phi(\theta^*, (x, \overline{\mathbf{Y}}_{i-1}), (x', Y_i))\| \rho^{i+m-1}.$$

Note that this term is small when $i$ is far from $-m$, say, $i \geq (k-m)/2$. Another kind of inequality is required to bound $\|\overline{\mathbb{E}}_{\theta^*}[\phi(\theta^*, \overline{\mathbf{Z}}_{i-1}^i)|\overline{\mathbf{Y}}_{-m}^k] - \overline{\mathbb{E}}_{\theta^*}[\phi(\theta^*, \overline{\mathbf{Z}}_{i-1}^i)|\overline{\mathbf{Y}}_{-m}^{k-1}]\|$. This type of bound will follow from forgetting properties of the reverse conditional hidden chain. Similar to Lemma 1, we have the following result.

LEMMA 9. *Assume* (A1) *and* (A2). *Let* $m, n \in \mathbb{Z}$ *with* $m, n \geq 0$ *and* $\theta \in \Theta$. *Under* $\overline{\mathbb{P}}_\theta$, *conditionally on* $\overline{\mathbf{Y}}_{-m}^n$, *the time-reversed process* $\{X_{n-k}\}_{0 \leq k \leq n+m}$ *is an inhomogeneous Markov chain, and for all* $0 < k \leq n+m$ *there exists a function* $\tilde{\mu}_k(\mathbf{y}_{-m-s+1}^{n-k}, A)$ *such that:*

(i) *for any* $A \in \mathcal{B}(\mathcal{X})$, $\mathbf{y}_{-m-s+1}^{n-k} \mapsto \tilde{\mu}_k(\mathbf{y}_{-m-s+1}^{n-k}, A)$ *is a Borel function;*

(ii) *for any* $\mathbf{y}_{-m-s+1}^{n-k}$, $\tilde{\mu}_k(\mathbf{y}_{-m-s+1}^{n-k}, \cdot)$ *is a probability measure on* $\mathcal{B}(\mathcal{X})$.
*In addition, for all* $\mathbf{y}_{-m-s+1}^{n-k}$, $\tilde{\mu}_k(\mathbf{y}_{-m-s+1}^{n-k}, \cdot) \ll \mu$ *and for all* $\overline{\mathbf{Y}}_{-m}^n$,

$$\overline{\mathbb{P}}_\theta(X_{n-k} \in A | X_{n-k+1}, \overline{\mathbf{Y}}_{-m}^n) = \overline{\mathbb{P}}_\theta(X_{n-k} \in A | X_{n-k+1}, \overline{\mathbf{Y}}_{-m}^{n-k})$$

$$\geq \frac{\sigma_-}{\sigma_+} \tilde{\mu}_k(\mathbf{Y}_{-m-s+1}^{n-k}, A).$$

The proof is along the same lines as Lemma 1 and is omitted for brevity.

From this lemma, using an analogue of Corollary 1, it follows that for $-m < i < k$,

$$\begin{aligned}(20)\quad &\|\overline{\mathbb{E}}_{\theta^*}[\phi(\theta^*, \overline{\mathbf{Z}}_{i-1}^i) | \overline{\mathbf{Y}}_{-m}^k] - \overline{\mathbb{E}}_{\theta^*}[\phi(\theta^*, \overline{\mathbf{Z}}_{i-1}^i) | \overline{\mathbf{Y}}_{-m}^{k-1}]\| \\ &\leq 2 \sup_{x,x' \in \mathcal{X}} \|\phi(\theta^*, (x, \overline{\mathbf{Y}}_{i-1}), (x', Y_i))\| \rho^{k-i-1}.\end{aligned}$$

By a standard martingale theory result [see, e.g., Shiryaev (1996), page 510], under assumption (A7) $\overline{\mathbb{E}}_{\theta^*}[\phi(\theta^*, \overline{\mathbf{Z}}_{i-1}^i) | \overline{\mathbf{Y}}_{-m}^k] \to \overline{\mathbb{E}}_{\theta^*}[\phi(\theta^*, \overline{\mathbf{Z}}_{i-1}^i) | \overline{\mathbf{Y}}_{-\infty}^k]$, $\overline{\mathbb{P}}_{\theta^*}$-a.s. as $m \to \infty$. Hence inequalities (19) and (20) hold true, $\overline{\mathbb{P}}_{\theta^*}$-a.s., when either $m$ or $m'$ is replaced by $\infty$. Using (20) with $m = \infty$ shows that

$$\sum_{i=-\infty}^{k-1} \|\overline{\mathbb{E}}_{\theta^*}[\phi(\theta^*, \overline{\mathbf{Z}}_{i-1}^i) | \overline{\mathbf{Y}}_{-\infty}^k] - \overline{\mathbb{E}}_{\theta^*}[\phi(\theta^*, \overline{\mathbf{Z}}_{i-1}^i) | \overline{\mathbf{Y}}_{-\infty}^{k-1}]\|$$

$$\leq \sum_{i=-\infty}^{k-1} 2 \sup_{x,x' \in \mathcal{X}} \|\phi(\theta^*, (x, \overline{\mathbf{Y}}_{i-1}), (x', Y_i))\| \rho^{k-i-1}.$$

Under (A7) the right-hand side is in $L^2(\overline{\mathbb{P}}_{\theta^*})$, and we may thus define

$$\Delta_{k,\infty}(\theta^*) \triangleq \overline{\mathbb{E}}_{\theta^*}[\phi(\theta^*, \overline{\mathbf{Z}}_{k-1}^k) | \overline{\mathbf{Y}}_{-\infty}^k]$$

$$+ \sum_{i=-\infty}^{k-1} (\overline{\mathbb{E}}_{\theta^*}[\phi(\theta^*, \overline{\mathbf{Z}}_{i-1}^i) | \overline{\mathbf{Y}}_{-\infty}^k] - \overline{\mathbb{E}}_{\theta^*}[\phi(\theta^*, \overline{\mathbf{Z}}_{i-1}^i) | \overline{\mathbf{Y}}_{-\infty}^{k-1}]).$$

In addition we have the following $L^2$-bound, showing that $\Delta_{k,m}(\theta^*)$ converges to $\Delta_{k,\infty}(\theta^*)$ in $L^2(\overline{\mathbb{P}}_{\theta^*})$ as $m \to \infty$.

LEMMA 10. *Assume* (A1), (A2) *and* (A6)–(A7). *Then, for all* $k \geq 1$ *and* $m \geq 0$,

$$(\overline{\mathbb{E}}_{\theta^*} \|\Delta_{k,m}(\theta^*) - \Delta_{k,\infty}(\theta^*)\|^2)^{1/2}$$

$$\leq 12 \left(\overline{\mathbb{E}}_{\theta^*} \left[\sup_{x,x' \in \mathcal{X}} \|\phi(\theta^*, (x, \overline{\mathbf{Y}}_0), (x', Y_1))\|^2\right]\right)^{1/2} \frac{\rho^{(k+m)/2-1}}{1-\rho}.$$

PROOF. Using (19) and (20) and the Minkowski inequality, we find that apart from the factor $(\overline{\mathbb{E}}_{\theta^*}[\sup_{x,x'\in\mathcal{X}}\|\phi(\theta^*,(x,\overline{\mathbf{Y}}_0),(x',Y_1))\|^2])^{1/2}$, $(\overline{\mathbb{E}}_{\theta^*}\|\Delta_{k,m}(\theta^*) - \Delta_{k,\infty}(\theta^*)\|^2)^{1/2}$ is bounded by

$$2\rho^{k+m-1} + 4\sum_{i=-m+1}^{k-1}(\rho^{k-i-1}\wedge\rho^{i+m-1}) + 2\sum_{i=-\infty}^{-m}\rho^{k-i-1}$$

$$\leq 2\rho^{k+m-1} + 4\sum_{-\infty<i\leq(k-m)/2}\rho^{k-i-1} + 4\sum_{(k-m)/2\leq i<\infty}\rho^{i+m-1} + 2\frac{\rho^{k+m-1}}{1-\rho}$$

$$\leq 12\frac{\rho^{(k+m)/2-1}}{1-\rho}.\qquad\square$$

Now define the filtration $\mathcal{F}$ by $\mathcal{F}_k = \sigma(\overline{\mathbf{Y}}_i; -\infty < i \leq k)$ for $k \in \mathbb{Z}$. By the conditional dominated convergence theorem,

$$\overline{\mathbb{E}}_{\theta^*}\left[\sum_{i=-\infty}^{k-1}(\overline{\mathbb{E}}_{\theta^*}[\phi(\theta^*,\overline{\mathbf{Z}}_{i-1}^i)|\overline{\mathbf{Y}}_{-\infty}^k] - \overline{\mathbb{E}}_{\theta^*}[\phi(\theta^*,\overline{\mathbf{Z}}_{i-1}^i)|\overline{\mathbf{Y}}_{-\infty}^{k-1}])\Big|\overline{\mathbf{Y}}_{-\infty}^{k-1}\right] = 0,$$

$$\overline{\mathbb{E}}_{\theta^*}[\phi(\theta^*,\overline{\mathbf{Z}}_{k-1}^k)|\overline{\mathbf{Y}}_{-\infty}^{k-1}] = \overline{\mathbb{E}}_{\theta^*}[\overline{\mathbb{E}}_{\theta^*}[\phi(\theta^*,\overline{\mathbf{Z}}_{k-1}^k)|\overline{\mathbf{Y}}_{-\infty}^{k-1},X_{k-1}]|\overline{\mathbf{Y}}_{-\infty}^{k-1}] = 0,$$

so that $\{\Delta_{k,\infty}(\theta^*)\}_{k=-\infty}^{\infty}$ is an $(\mathcal{F},\overline{\mathbb{P}}_{\theta^*})$-adapted stationary, ergodic and square integrable martingale increment sequence. The CLT for sums of such sequences [see, e.g., Durrett (1996), page 418] shows that

$$n^{-1/2}\sum_{k=1}^{n}\Delta_{k,\infty}(\theta^*) \to \mathcal{N}(0,I(\theta^*)),\qquad \overline{\mathbb{P}}_{\theta^*}\text{-weakly},$$

where $I(\theta^*) \triangleq \overline{\mathbb{E}}_{\theta^*}[\Delta_{0,\infty}(\theta^*)\Delta_{0,\infty}(\theta^*)^T]$ is the asymptotic Fisher information matrix, defined as the covariance matrix of the asymptotic score function. Lemma 10 implies that

$$(21)\qquad \lim_{n\to\infty}\overline{\mathbb{E}}_{\theta^*}\left\|n^{-1/2}\sum_{k=1}^{n}(\Delta_{k,0}(\theta^*) - \Delta_{k,\infty}(\theta^*))\right\|^2 = 0,$$

and hence $n^{-1/2}\sum_{k=1}^n \Delta_{k,0}(\theta^*)$, and by Lemma 8 also $n^{-1/2}\sum_{k=1}^n \Delta_{k,0,x}(\theta^*)$, have the same limiting distribution under $\overline{\mathbb{P}}_{\theta^*}$. We summarize our findings in the following result.

THEOREM 2. *Assume* (A1), (A2) *and* (A6)–(A8). *Then for any* $x_0 \in \mathcal{X}$,

$$n^{-1/2}\nabla_\theta l_n(\theta^*,x_0) \to \mathcal{N}(0,I(\theta^*)),\qquad \overline{\mathbb{P}}_{\theta^*}\text{-weakly.}$$

6.2. *Law of large numbers for the observed Fisher information.* The second part of the proof consists of showing a locally uniform law of large numbers for the observed Fisher information; for all possibly random sequences $\{\theta_n^*\}$ such that $\theta_n^* \to \theta^*$, $\mathbb{P}_{\theta^*}$-a.s., $-n^{-1}\nabla_\theta^2 l_n(\theta_n^*, x_0)$ converges, $\mathbb{P}_{\theta^*}$-a.s., to the Fisher information matrix at $\theta^*$. Similar to what was done in the previous section and following the ideas developed in BRR, the proof amounts to showing that $-n^{-1}\nabla_\theta^2 l_n(\theta_n^*, x_0)$ may be approximated by the sample mean of an ergodic stationary process. To do that it is convenient, just as for the score function, to express the observed Fisher information in terms of the Hessian of the complete log likelihood. This can be done by using the so-called Louis missing information principle [Louis (1982)], valid under (A7) and (A8), which shows that

$$
\begin{aligned}
\nabla_\theta^2 \log p_\theta(\mathbf{Y}_1^n | \overline{\mathbf{Y}}_0, X_0 = x_0) \\
= \mathbb{E}_\theta \left[ \sum_{i=1}^n \varphi(\theta, \overline{\mathbf{Z}}_{i-1}^i) \Big| \overline{\mathbf{Y}}_0^n, X_0 = x_0 \right] \\
+ \operatorname{var}_\theta \left[ \sum_{i=1}^n \phi(\theta, \overline{\mathbf{Z}}_{i-1}^i) \Big| \overline{\mathbf{Y}}_0^n, X_0 = x_0 \right],
\end{aligned}
\tag{22}
$$

where

$$
\begin{aligned}
\varphi(\theta, \overline{\mathbf{Z}}_{i-1}^i) &= \varphi(\theta, \overline{\mathbf{Z}}_{i-1}, \overline{\mathbf{Z}}_i) = \varphi(\theta, (\mathbf{X}_{i-s}^{i-1}, \overline{\mathbf{Y}}_{i-1}), (\mathbf{X}_{i-s+1}^i, \mathbf{Y}_{i-s+1}^i)) \\
&= \varphi(\theta, (X_{i-1}, \overline{\mathbf{Y}}_{i-1}), (X_i, Y_i)) \\
&\triangleq \nabla_\theta^2 \log\left(q_\theta(X_{i-1}, X_i) g_\theta(Y_i | \overline{\mathbf{Y}}_{i-1}, X_i)\right).
\end{aligned}
$$

As above we may write these quantities as telescoping sums:

$$
\begin{aligned}
&\mathbb{E}_\theta \left[ \sum_{i=1}^n \varphi(\theta, \overline{\mathbf{Z}}_{i-1}^i) \Big| \overline{\mathbf{Y}}_0^n, X_0 = x_0 \right] \\
&= \sum_{k=1}^n \left( \mathbb{E}_\theta \left[ \sum_{i=1}^k \varphi(\theta, \overline{\mathbf{Z}}_{i-1}^i) \Big| \overline{\mathbf{Y}}_0^k, X_0 = x_0 \right] \right. \\
&\quad \left. - \mathbb{E}_\theta \left[ \sum_{i=1}^{k-1} \varphi(\theta, \overline{\mathbf{Z}}_{i-1}^i) \Big| \overline{\mathbf{Y}}_0^{k-1}, X_0 = x_0 \right] \right)
\end{aligned}
$$

and

$$\operatorname{var}_\theta \left[ \sum_{i=1}^n \phi(\theta, \overline{\mathbf{Z}}_{i-1}^i) \Big| \overline{\mathbf{Y}}_0^n, X_0 = x_0 \right]$$

$$= \sum_{k=1}^n \left( \operatorname{var}_\theta \left[ \sum_{i=1}^k \phi(\theta, \overline{\mathbf{Z}}_{i-1}^i) \Big| \overline{\mathbf{Y}}_0^k, X_0 = x_0 \right] \right.$$

$$\left. - \operatorname{var}_\theta \left[ \sum_{i=1}^{k-1} \phi(\theta, \overline{\mathbf{Z}}_{i-1}^i) \Big| \overline{\mathbf{Y}}_0^{k-1}, X_0 = x_0 \right] \right).$$

It turns out (see Lemma 13) that as $k \to \infty$ the initial condition on $X_0$ becomes irrelevant. Therefore it is sensible to define, for $k \geq 1$ and $m \geq 0$,

$$\Delta_{k,m}(\theta) = \overline{\mathbb{E}}_\theta \left[ \sum_{i=-m+1}^k \varphi(\theta, \overline{\mathbf{Z}}_{i-1}^i) \Big| \overline{\mathbf{Y}}_{-m}^k \right]$$
(23)
$$- \overline{\mathbb{E}}_\theta \left[ \sum_{i=-m+1}^{k-1} \varphi(\theta, \overline{\mathbf{Z}}_{i-1}^i) \Big| \overline{\mathbf{Y}}_{-m}^{k-1} \right],$$

$$\Gamma_{k,m}(\theta) = \overline{\operatorname{var}}_\theta \left[ \sum_{i=-m+1}^k \phi(\theta, \overline{\mathbf{Z}}_{i-1}^i) \Big| \overline{\mathbf{Y}}_{-m}^k \right]$$
(24)
$$- \overline{\operatorname{var}}_\theta \left[ \sum_{i=-m+1}^{k-1} \phi(\theta, \overline{\mathbf{Z}}_{i-1}^i) \Big| \overline{\mathbf{Y}}_{-m}^{k-1} \right].$$

Propositions 4 and 5 show that $\Delta_{k,m}(\theta)$ and $\Gamma_{k,m}(\theta)$ both have limits as $m \to \infty$, $\overline{\mathbb{P}}_{\theta^*}$-a.s., and in $L^1(\overline{\mathbb{P}}_{\theta^*})$. Let $\Delta_{k,\infty}(\theta)$ and $\Gamma_{k,\infty}(\theta)$ denote these limits. It follows from the definitions above that $\{\Delta_{k,\infty}\}_{k=1}^\infty$ and $\{\Gamma_{k,\infty}\}_{k=1}^\infty$ are $\overline{\mathbb{P}}_{\theta^*}$-stationary and ergodic, and the limit of the observed Fisher information will be $-\overline{\mathbb{E}}_{\theta^*}[\Delta_{0,\infty}(\theta^*) + \Gamma_{0,\infty}(\theta^*)]$.

PROPOSITION 4. *Assume* (A1)–(A3). *Let $G$ be a compact subset of $\Theta$, let $q > 0$ and let $\varphi \colon \Theta \times \mathcal{X}^q \times \mathcal{Y}^q \to \mathbb{R}$ be a Borel function such that for all $\mathbf{x}_1^q \in \mathcal{X}^q$ and $\mathbf{y}_1^q \in \mathcal{Y}^q$, $\varphi(\theta, \mathbf{x}_1^q, \mathbf{y}_1^q)$ is continuous w.r.t. $\theta$ on $G$ and*

$$\overline{\mathbb{E}}_{\theta^*} \left[ \sup_{\theta \in G} \sup_{\mathbf{x}_1^q \in \mathcal{X}^q} |\varphi(\theta, \mathbf{x}_1^q, \mathbf{Y}_1^q)| \right] < \infty.$$

*Then for each $\theta \in G$, $\Delta_{k,m}(\theta)$, as defined in (23), converges $\overline{\mathbb{P}}_{\theta^*}$-a.s. and in $L^1(\overline{\mathbb{P}}_{\theta^*})$ to $\Delta_{k,\infty}(\theta)$ as $m \to \infty$. In addition, the function $\theta \mapsto \overline{\mathbb{E}}_{\theta^*}[\Delta_{0,\infty}(\theta)]$ is continuous on $G$ and for all $x_0 \in \mathcal{X}$ and $\theta \in G$,*

$$\lim_{\delta \to 0} \lim_{n \to \infty} \sup_{|\theta' - \theta| \leq \delta} \left| n^{-1} \mathbb{E}_{\theta'} \left[ \sum_{i=1}^n \varphi(\theta', \mathbf{Z}_{i-q+1}^i) \Big| \overline{\mathbf{Y}}_0^n, X_0 = x_0 \right] \right.$$

$$\left. - \overline{\mathbb{E}}_{\theta^*}[\Delta_{0,\infty}(\theta)] \right| = 0, \qquad \overline{\mathbb{P}}_{\theta^*}\text{-}a.s.$$

PROPOSITION 5. *Assume* (A1)–(A3). *Let $G$ be a compact subset of $\Theta$, let $q > 0$ and let $\phi : \Theta \times \mathcal{X}^q \times \mathcal{Y}^q \to \mathbb{R}$ be a Borel function such that for all $\mathbf{x}_1^q \in \mathcal{X}^q$ and $\mathbf{y}_1^q \in \mathcal{Y}^q$, $\phi(\theta, \mathbf{x}_1^q, \mathbf{y}_1^q)$ is continuous w.r.t. $\theta$ on $G$ and*

$$\overline{\mathbb{E}}_{\theta^*}\left[\sup_{\theta \in G} \sup_{\mathbf{x}_1^q \in \mathcal{X}^q} |\phi(\theta, \mathbf{x}_1^q, \mathbf{Y}_1^q)|^2\right] < \infty.$$

*Then for each $\theta \in G$, $\Gamma_{k,m}(\theta)$, as defined in* (24), *converges $\overline{\mathbb{P}}_{\theta^*}$-a.s. and in $L^1(\overline{\mathbb{P}}_{\theta^*})$ to $\Gamma_{k,\infty}(\theta)$ as $m \to \infty$. In addition, the function $\theta \mapsto \overline{\mathbb{E}}_{\theta^*}[\Gamma_{0,\infty}(\theta)]$ is continuous on $G$ and, for all $x_0 \in \mathcal{X}$ and $\theta \in G$,*

$$\lim_{\delta \to 0} \lim_{n \to \infty} \sup_{|\theta' - \theta| \leq \delta} \left| n^{-1} \mathrm{var}_{\theta'}\left[\sum_{i=1}^n \phi(\theta', \mathbf{Z}_{i-q+1}^i) \middle| \overline{\mathbf{Y}}_0^n, X_0 = x_0\right] \right.$$
$$\left. - \overline{\mathbb{E}}_{\theta^*}[\Gamma_{0,\infty}(\theta)] \right| = 0, \qquad \overline{\mathbb{P}}_{\theta^*}\text{-}a.s.$$

Note that in Propositions 4 and 5 the functions $\varphi$ and $\phi$ take values in $\mathbb{R}$. Adaptations to vector- and matrix-valued functions are straightforward.

For all $x_0 \in \mathcal{X}$ the Fisher information identity implies, under the stated assumptions, that

$$n^{-1}\mathbb{E}_\theta[\nabla_\theta l_n(\theta, x_0) \nabla_\theta l_n(\theta, x_0)^T | \overline{\mathbf{Y}}_0, X_0 = x_0]$$
$$= -n^{-1}\mathbb{E}_\theta[\nabla_\theta^2 l_n(\theta, x_0) | \overline{\mathbf{Y}}_0, X_0 = x_0],$$

and Propositions 4 and 5 together with the Louis missing information principle show that the limits in $n$ of these two quantities both coincide with the Fisher information at $\theta^*$. We conclude the discussion in this section by stating the main result.

THEOREM 3. *Assume* (A1)–(A3) *and* (A6)–(A8) *and let $\{\theta_n^*\}$ be any, possibly stochastic, sequence in $\Theta$ such that $\theta_n^* \to \theta^*$ $\overline{\mathbb{P}}_{\theta^*}$-a.s. Then for all $x_0 \in \mathcal{X}$ $-n^{-1}\nabla_\theta^2 l_n(\theta_n^*, x_0) \to I(\theta^*)$, $\overline{\mathbb{P}}_{\theta^*}$-a.s.*

The following theorem is a standard consequence of Theorems 2 and 3 (see, e.g., BRR).

THEOREM 4. *Assume* (A1)–(A8) *and that $I(\theta^*)$ is positive definite. Then for all $x_0 \in \mathcal{X}$*

$$n^{1/2}(\hat{\theta}_{n,x_0} - \theta^*) \to \mathcal{N}(0, I(\theta^*)^{-1}), \qquad \overline{\mathbb{P}}_{\theta^*}\text{-}weakly.$$

**7. Extensions to nonstationary AR models with Markov regime.** In Sections 4 and 6 the assumption of stationarity of $\{Y_k\}$ plays a crucial role. In this section we shall extend the consistency and asymptotic normality of the MLE to the case where this process is not stationary. Hence we assume that the process we observe, denoted by $\{\overline{\mathbf{Y}}'_k\}_{k=0}^\infty$, and the associated hidden chain, denoted by $\{X'_k\}_{k=0}^\infty$, are governed by the transition kernel $\Pi_{\theta^*}$ and with $(X'_0, \overline{\mathbf{Y}}'_0)$ having distribution $\zeta$. This initial distribution is unknown to us and in general $\zeta \neq \pi_{\theta^*}$. As before we let $\{(X_k, \overline{\mathbf{Y}}_k)\}_{k=0}^\infty$ denote a corresponding stationary process.

We observe that since these processes are positive Harris recurrent and aperiodic [this is (A2)] we can construct them on a common probability space in a way that there exists an a.s. finite random time $T$, the coupling time, such that $\overline{\mathbf{Z}}_n = \overline{\mathbf{Z}}'_n$ for $n \geq T$ [Thorisson (2000), page 369]. The associated probability measure is denoted by $\mathbb{P}_{\pi_{\theta^*} \otimes \zeta}$. Hence, to be precise, $\mathbb{P}_{\pi_{\theta^*} \otimes \zeta}(T < \infty) = 1$.

Define $l'_n(\theta, x_0) \triangleq \log p_\theta([\mathbf{Y}']_1^n | \overline{\mathbf{Y}}'_0, X'_0 = x_0)$ and let $\hat{\theta}'_{n,x_0}$ be the maximizer of this function w.r.t. $\theta$. Put

$$\begin{aligned} D_n(\theta, x_0) &\triangleq l'_n(\theta, x_0) - l_n(\theta, x_0) \\ &= \sum_{k=1}^n (\log p_\theta(Y'_k | [\overline{\mathbf{Y}}']_0^{k-1}, X'_0 = x_0) - \log p_\theta(Y_k | \overline{\mathbf{Y}}_0^{k-1}, X_0 = x_0)). \end{aligned}$$

The following lemma ensures that $D_n(\theta, x_0)$ is bounded, $\mathbb{P}_{\pi_{\theta^*} \otimes \zeta}$-a.s., which implies that the difference between $\hat{\theta}'_{n,x_0}$ and $\hat{\theta}_{n,x_0}$ converges to zero, $\mathbb{P}_{\pi_{\theta^*} \otimes \zeta}$-a.s. (see Theorem 5).

LEMMA 11. *Assume* (A1) *and* (A2). *Then for all* $\zeta$ *and all* $x_0 \in \mathcal{X}$, $\sup_{n \geq 0} \sup_{\theta \in \Theta} |D_n(\theta, x_0)| < \infty$, $\mathbb{P}_{\pi_{\theta^*} \otimes \zeta}$-*a.s.*

PROOF. Write

$$\sup_{\theta \in \Theta} |D_n(\theta, x_0)|$$

$$\leq \sum_{k=1}^\infty \sup_{\theta \in \Theta} |\log p_\theta(Y'_k | [\overline{\mathbf{Y}}']_0^{k-1}, X'_0 = x_0) - \log p_\theta(Y_k | \overline{\mathbf{Y}}_0^{k-1}, X_0 = x_0)|$$

$$\leq \sum_{k=1}^T \left( \sup_{\theta \in \Theta} |\log p_\theta(Y'_k | [\overline{\mathbf{Y}}']_0^{k-1}, X'_0 = x_0)| \right.$$

(25) $$\left. + \sup_{\theta \in \Theta} |\log p_\theta(Y_k | \overline{\mathbf{Y}}_0^{k-1}, X_0 = x_0)| \right)$$

$$+ \sum_{k=T+1}^{\infty} \sup_{\theta \in \Theta} |\log p_\theta(Y'_k|[\overline{\mathbf{Y}}']_0^{k-1}, X'_0 = x_0)$$

$$- \log p_\theta(Y_k|\overline{\mathbf{Y}}_0^{k-1}, X_0 = x_0)|.$$

Since

$$\sigma_- \int g_\theta(Y_k|\overline{\mathbf{Y}}_{k-1}, x)\mu(dx) \le p_\theta(Y_k|\overline{\mathbf{Y}}_0^{k-1}, X_0 = x_0)$$

$$\le \sigma_+ \int g_\theta(Y_k|\overline{\mathbf{Y}}_{k-1}, x)\mu(dx)$$

(see the proof of Lemma 2), the first sum on the right-hand side is finite $\mathbb{P}_{\pi_{\theta^*} \otimes \zeta}$-a.s. by (A1).

For the second sum, note that for all $i < k$,

$$p_\theta(Y_k|\overline{\mathbf{Y}}_0^{k-1}, X_0 = x_0)$$
$$= \int\int\int g_\theta(Y_k|\overline{\mathbf{Y}}_{k-1}, x_k) q_\theta(x_{k-1}, x_k) \mu(dx_k) \mathbb{P}_\theta(dx_{k-1}|x_i, \overline{\mathbf{Y}}_i^{k-1})$$
$$\times \mathbb{P}_\theta(dx_i|\overline{\mathbf{Y}}_0^{k-1}, X_0 = x_0),$$

and similarly for $p_\theta(Y'_k|[\overline{\mathbf{Y}}']_0^{k-1}, X'_0 = x_0)$. Using the fact that for $n \ge T$, $\overline{\mathbf{Z}}_n = \overline{\mathbf{Z}}'_n$ and thus $\overline{\mathbf{Y}}_n = \overline{\mathbf{Y}}'_n$ and Corollary 1, we have for all $k > T$,

$$|p_\theta(Y_k|[\overline{\mathbf{Y}}']_0^{k-1}, X'_0 = x_0) - p_\theta(Y_k|\overline{\mathbf{Y}}_0^{k-1}, X_0 = x_0)|$$
$$\le \rho^{k-T-1} \sigma_+ \int g_\theta(Y_k|\overline{\mathbf{Y}}_{k-1}, x)\mu(dx),$$

and hence

(26)
$$|\log p_\theta(Y'_k|[\overline{\mathbf{Y}}']_0^{k-1}, X'_0 = x_0) - \log p_\theta(Y_k|\overline{\mathbf{Y}}_0^{k-1}, X_0 = x_0)|$$
$$\le \rho^{k-T-1}/(1-\rho);$$

compare the proof of Lemma 2. Thus the second sum on the right-hand side of (25) is also finite $\mathbb{P}_{\pi_{\theta^*} \otimes \zeta}$-a.s. □

We now can prove the consistency of the MLE for a nonstationary process.

THEOREM 5. *Assume* (A1)–(A5). *Then for all $\zeta$ and any $x_0 \in \mathcal{X}$, $\lim_{n\to\infty} \hat{\theta}'_{n,x_0} = \theta^*$, $\mathbb{P}_{\pi_{\theta^*} \otimes \zeta}$-a.s.*

PROOF. Since $\hat{\theta}'_{n,x_0}$ is the maximizer of $\theta \mapsto n^{-1}l'_n(\theta,x_0)$,

$$\begin{aligned}
l(\hat{\theta}'_{n,x_0}) &\geq l(\theta^*) - l(\theta^*) + n^{-1}l_n(\theta^*,x_0) \\
&\quad - n^{-1}l_n(\theta^*,x_0) + n^{-1}l'_n(\theta^*,x_0) - n^{-1}l'_n(\hat{\theta}'_{n,x_0},x_0) \\
&\quad + n^{-1}l_n(\hat{\theta}'_{n,x_0},x_0) - n^{-1}l_n(\hat{\theta}'_{n,x_0},x_0) + l(\hat{\theta}'_{n,x_0}) \\
&\geq l(\theta^*) - 2\sup_{\theta\in\Theta}|n^{-1}l_n(\theta,x_0) - l(\theta)| - 2\sup_{\theta\in\Theta}|n^{-1}D_n(\theta,x_0)|.
\end{aligned}$$

The right-hand side of this inequality tends to $l(\theta^*)$, $\mathbb{P}_{\pi_{\theta^*}\otimes\zeta}$-a.s., by Proposition 2 and Lemma 11. The proof now follows from Proposition 3, continuity of $l(\theta)$ (Lemma 4) and compactness of $\Theta$. □

To show that $n^{1/2}(\hat{\theta}'_{n,x_0} - \hat{\theta}_{n,x_0}) \to 0$, $\mathbb{P}_{\pi_{\theta^*}\otimes\zeta}$-a.s. and thus that $\hat{\theta}'_{n,x_0}$ and $\hat{\theta}_{n,x_0}$ are asymptotically normal with the same covariance matrix, we need to show some kind of continuity of the function $\theta \mapsto D_n(\theta,x_0)$.

LEMMA 12. *Assume* (A1)–(A5). *Then*

$$\lim_{n\to\infty}|D_n(\hat{\theta}'_{n,x_0},x_0) - D_n(\hat{\theta}_{n,x_0},x_0)| = 0, \qquad \mathbb{P}_{\pi_{\theta^*}\otimes\zeta}\text{-}a.s.$$

PROOF. Let $\varepsilon > 0$. By (26) there exists a random integer $N$ which is finite $\mathbb{P}_{\pi_{\theta^*}\otimes\zeta}$-a.s. and satisfies

$$\sum_{k=N+1}^{\infty}\sup_{\theta\in\Theta}|\log p_\theta(Y'_k|[\overline{\mathbf{Y}}']_0^{k-1}) - \log p_\theta(Y_k|\overline{\mathbf{Y}}_0^{k-1})| \leq \varepsilon, \qquad \mathbb{P}_{\pi_{\theta^*}\otimes\zeta}\text{-a.s.}$$

Thus, $\mathbb{P}_{\pi_{\theta^*}\otimes\zeta}$-a.s. for all $n \geq N$,

$$\begin{aligned}
&|D_n(\hat{\theta}'_{n,x_0},x_0) - D_n(\hat{\theta}_{n,x_0},x_0)| \\
&\quad \leq 2\varepsilon + |l'_N(\hat{\theta}'_{n,x_0},x_0) - l'_N(\hat{\theta}_{n,x_0},x_0)| + |l_N(\hat{\theta}'_{n,x_0},x_0) - l_N(\hat{\theta}_{n,x_0},x_0)|.
\end{aligned}$$

Under the given assumptions $\theta \mapsto l'_N(\theta,x_0)$ and $\theta \mapsto l_N(\theta,x_0)$ are $\mathbb{P}_{\pi_{\theta^*}\otimes\zeta}$-a.s. continuous (see the proof of Lemma 4) and the proof is complete upon observing that $\hat{\theta}'_{n,x_0}$ and $\hat{\theta}_{n,x_0}$ both converge to $\theta^*$, $\mathbb{P}_{\pi_{\theta^*}\otimes\zeta}$-a.s., and that $\varepsilon$ was arbitrary. □

THEOREM 6. *Assume* (A1)–(A8) *and that* $I(\theta^*)$ *is positive definite. Then, for all* $\zeta$ *and any* $x_0 \in \mathcal{X}$,

$$n^{1/2}(\hat{\theta}'_{n,x_0} - \theta^*) \to \mathcal{N}(0,I(\theta^*)^{-1}), \qquad \mathbb{P}_{\pi_{\theta^*}\otimes\zeta}\text{-weakly}.$$

PROOF. It is sufficient to prove that $\varepsilon_n \triangleq \sqrt{n}(\hat{\theta}_{n,x_0} - \hat{\theta}'_{n,x_0}) \to 0$, $\mathbb{P}_{\pi_{\theta^*} \otimes \zeta}$-a.s. Since $\hat{\theta}'_{n,x_0}$ is the maximizer of $\theta \mapsto l'_n(\theta, x_0)$, $l'_n(\hat{\theta}'_{n,x_0}, x_0) \geq l'_n(\hat{\theta}_{n,x_0}, x_0)$, which implies that

$$D_n(\hat{\theta}'_{n,x_0}, x_0) - D_n(\hat{\theta}_{n,x_0}, x_0) \geq l_n(\hat{\theta}_{n,x_0}, x_0) - l_n(\hat{\theta}'_{n,x_0}, x_0)$$
$$= -\tfrac{1}{2} n^{-1} \varepsilon_n^T \nabla_\theta^2 l_n(t_n \hat{\theta}'_{n,x_0} + (1-t_n)\hat{\theta}_{n,x_0}) \varepsilon_n$$

for some $0 \leq t_n \leq 1$. By a straightforward adaptation of Theorem 3 to the present case with two processes,

$$-n^{-1} \nabla_\theta^2 l_n(t_n \hat{\theta}'_{n,x_0} + (1-t_n) \hat{\theta}_{n,x_0}) \to I(\theta^*), \qquad \mathbb{P}_{\pi_{\theta^*} \otimes \zeta}\text{-a.s.}$$

Since $I(\theta^*)$ is positive definite there exists $M > 0$ such that on a set with $\mathbb{P}_{\pi_{\theta^*} \otimes \zeta}$-probability one and for $n$ sufficiently large,

$$D_n(\hat{\theta}'_{n,x_0}, x_0) - D_n(\hat{\theta}_{n,x_0}, x_0) \geq M |\varepsilon_n|^2.$$

The proof is complete by applying Lemma 12. □

## 8. Numerical approximations.

8.1. *Two Monte Carlo numerical methods.* As mentioned in the Introduction, when the state space of $\{X_k\}$ is continuous the log likelihood needs to be approximated by some numerical method. Here we list two classes of methods that have been proposed and successfully used in many practical problems, but point out that there are other ones as well, for example, importance sampling [Geyer and Thompson (1992) and Geyer (1994)].

*Particle filters.* These methods depart from the representation

$$l_n(\theta, x_0) = \sum_{k=1}^n \log \int g_\theta(Y_k | \overline{\mathbf{Y}}_{k-1}, x_k) \mathbb{P}_\theta(X_k \in dx_k | \overline{\mathbf{Y}}_0^{k-1}, X_0 = x_0)$$

and replace the predictive distribution $\mathbb{P}_\theta(X_k \in dx_k | \overline{\mathbf{Y}}_0^{k-1}, X_0 = x_0)$ by a particle approximation. More precisely, the approximating distribution is the empirical distribution of the locations of $N$ particles at time $k$. There are many variants to how the locations of the particles are updated, and under general assumptions the particle approximation converges to the true predictive distribution at rate $N^{-1/2}$ when $N$ grows. The approximate log likelihood may be maximized using any standard numerical optimization algorithm. Further reading is found in the collection Doucet, de Freitas and Gordon (2001); see in particular the survey paper Hürzeler and Künsch (2001). Other references are Künsch (2001) and Pitt (2002). Particle filter methods have been proved to perform well in a wide range of problems, as illustrated in the above references.

*Monte Carlo EM algorithms.* The EM algorithm is an iterative algorithm for computing the MLE (or at least a local maximum of the log likelihood) in problems with missing data. Its key components are the computation of the function $Q(\theta, \theta') = \mathbb{E}_\theta[\log p_{\theta'}(\mathbf{X}_1^n, \mathbf{Y}_1^n | \overline{\mathbf{Y}}_0^n, X_0 = x_0) | \overline{\mathbf{Y}}_0^n, X_0 = x_0]$ (the E-step) and the maximization of this function w.r.t. $\theta'$ (the M-step). These two steps constitute the update from a current estimate $\theta$ to a new one. Obviously the EM user is required to compute conditional expectations of functions of $\mathbf{X}_1^n$ given $\overline{\mathbf{Y}}_0^n$ and $X_0 = x_0$. If the state space is continuous this task is typically infeasible, but the conditional expectations can be replaced by sample averages over $m$ simulated realizations of $\mathbf{X}_1^n$ under the same conditions. These methods are called Monte Carlo EM (MCEM) algorithms, or stochastic EM (SEM) algorithms. A recent survey is found in Booth, Hobert and Jank (2001), and general versions of the algorithm are described in Tanner (1996) and Nielsen (2000). If the number $m$ of simulated replications is allowed to increase with each iteration, the algorithm can be made to converge [Fort and Moulines (2003)]. MCEM methods are successfully used in many areas; see the above-mentioned survey paper.

Having said that, we stress that the distinction between particle filter and MCEM methods is not sharp. In fact, the function $Q(\theta, \theta')$ of the EM algorithm can, in principle, be computed recursively in $n$, which opens up for particle approximations of this functional [Cappé (2001)]. Hence, the approximation and maximization of the log likelihood rather splits into two other subproblems to be considered. First, the optimization scheme: (i) EM type, which is particularly appropriate if the complete data is from an exponential or curved exponential family of distributions, or (ii) a standard numerical optimization algorithm such as a quasi-Newton or conjugate gradient method. Second, the approach to approximate conditional expectations: (i) forward in time using particle filters or (ii) conditional on the whole set of data using more traditional MCMC simulation.

8.2. *Asymptotics of approximate estimators.* Theorems 1 and 4 give the asymptotic properties of the MLE, but, as noted above, neither the (conditional) likelihood nor the MLE is computable unless the state space is finite. An important question is thus if an approximate computation of the MLE or likelihood is sufficient to retain the asymptotics. Of course, if $\tilde{\theta}_{n,x_0}$ is an estimator such that $\tilde{\theta}_{n,x_0} - \hat{\theta}_{n,x_0} = o_P(n^{-1/2})$ (with $P = \overline{\mathbb{P}}_{\theta^*}$), then $\tilde{\theta}_{n,x_0}$ is consistent and $n^{1/2}(\tilde{\theta}_{n,x_0} - \theta^*)$ has the same distributional limit as $n^{1/2}(\hat{\theta}_{n,x_0} - \theta^*)$. This simple observation applies to methods that directly approximate the MLE, for example, MCEM. The following theorem gives a corresponding result when the likelihood is approximated.

THEOREM 7. *Assume that $\tilde{\theta}_{n,x_0}$ is an estimator satisfying $l_n(\tilde{\theta}_{n,x_0}, x_0) \geq \sup_{\theta \in \Theta} l_n(\theta, x_0) - R_n$ and that the assumptions of Theorem 4 hold. Then the following are true*:

(i) *If $R_n = o_P(n)$ (with $P = \overline{\mathbb{P}}_{\theta^*}$), then $\tilde{\theta}_{n,x_0}$ is consistent.*

(ii) *If $R_n = O_P(1)$, then $n^{1/2}(\tilde{\theta}_{n,x_0} - \theta^*) = O_P(1)$, that is, the sequence $\{\tilde{\theta}_{n,x_0}\}$ is $n^{1/2}$-consistent under $\overline{\mathbb{P}}_{\theta^*}$.*

(iii) *If $R_n = o_P(1)$, then $n^{1/2}(\tilde{\theta}_{n,x_0} - \theta^*) \to \mathcal{N}(0, I(\theta^*)^{-1})$, $\overline{\mathbb{P}}_{\theta^*}$-weakly as $n \to \infty$.*

REMARK 7. The remainder term does not depend on $\theta$, that is, it is uniform in $\theta \in \Theta$. If $n^{-1}R_n \to 0$, $\overline{\mathbb{P}}_{\theta^*}$-a.s. in (i), we obtain strong consistency.

PROOF OF THEOREM 7. We start with (i). Since $l_n(\tilde{\theta}_{n,x_0}, x_0) \geq \sup_{\theta \in \Theta} l_n(\theta, x_0) - R_n \geq l_n(\theta^*, x_0) - R_n$, we have

$$\begin{aligned} l(\theta^*) &\geq l(\tilde{\theta}_{n,x_0}) \\ &\geq l(\theta^*) - l(\theta^*) + n^{-1}l_n(\theta^*, x_0) - n^{-1}l_n(\tilde{\theta}_{n,x_0}, x_0) + l(\tilde{\theta}_{n,x_0}) - n^{-1}R_n \\ &\geq l(\theta^*) - 2\sup_{\theta \in \Theta}|n^{-1}l_n(\theta, x_0) - l(\theta)| - n^{-1}R_n. \end{aligned}$$

If $R_n = o_P(n)$, using Proposition 2, $l(\tilde{\theta}_{n,x_0}) - l(\theta^*) = o_P(1)$. Standard compactness arguments going back to Wald (1949) and Proposition 3 complete the proof of (i).

We now turn to (ii) and (iii). Recall that $\hat{\theta}_{n,x_0}$ maximizes $l_n(\theta, x_0)$. By a Taylor expansion of $l_n(\theta, x_0)$ around $\hat{\theta}_{n,x_0}$, there exists a point $\bar{\theta}_n$ on the line segment between $\hat{\theta}_{n,x_0}$ and $\tilde{\theta}_{n,x_0}$ such that

$$R_n \geq l_n(\hat{\theta}_{n,x_0}, x_0) - l_n(\tilde{\theta}_{n,x_0}, x_0) = \varepsilon_n^T(-n^{-1}\nabla_\theta^2 l_n(\bar{\theta}_n, x_0))\varepsilon_n,$$

where $\varepsilon_n = n^{1/2}(\tilde{\theta}_{n,x_0} - \hat{\theta}_{n,x_0})$. Since $\tilde{\theta}_{n,x_0}$ converges to $\theta^*$ in probability, so does $\bar{\theta}_n$. Hence there is a positive sequence $\{\delta_n\}$ tending to zero such that $\overline{\mathbb{P}}_{\theta^*}(|\bar{\theta}_n - \theta^*| > \delta_n) \to 0$. Thus, for any $c > 0$,

$$\begin{aligned} &\overline{\mathbb{P}}_{\theta^*}(\|-n^{-1}\nabla_\theta^2 l_n(\bar{\theta}_n, x_0) - I(\theta^*)\| > c) \\ &\quad \leq \overline{\mathbb{P}}_{\theta^*}(|\bar{\theta}_n - \theta^*| > \delta_n) + \overline{\mathbb{P}}_{\theta^*}\bigg(\sup_{|\theta - \theta^*| \leq \delta_n} \|-n^{-1}\nabla_\theta^2 l_n(\theta, x_0) - I(\theta^*)\| > c\bigg). \end{aligned}$$

The first term on the right-hand side tends to zero as $n \to \infty$, and so does the second one by Theorem 3. Since $I(\theta^*)$ is assumed positive definite, there exists an $M > 0$ such that

$$R_n \geq (M + o_P(1))|\varepsilon_n|^2.$$

Thus, if $R_n = O_P(1)$, then $\varepsilon_n = O_P(1)$, and if $R_n = o_P(1)$, then $\varepsilon_n = o_P(1)$. The proofs of (ii) and (iii) are now complete using $n^{1/2}(\tilde{\theta}_{n,x_0} - \theta^*) = \varepsilon_n + n^{1/2}(\hat{\theta}_{n,x_0} - \theta^*)$ and the result of Theorem 4. □

Obviously, if $\tilde{l}_n$ is an approximation of the true log likelihood $l_n$ (both conditional on $x_0$) such that $|\tilde{l}_n(\theta, x_0) - l_n(\theta, x_0)| \leq (1/2) R_n$ for all $\theta \in \Theta$, and $\tilde{\theta}_{n,x_0}$ is the corresponding maximizer, then $l_n(\tilde{\theta}_{n,x_0}, x_0) \geq \tilde{l}_n(\tilde{\theta}_{n,x_0}, x_0) - (1/2) R_n \geq \tilde{l}_n(\hat{\theta}_{n,x_0}, x_0) - (1/2) R_n \geq l_n(\hat{\theta}_{n,x_0}, x_0) - R_n$, that is, the principal condition of the theorem is fulfilled. We thus see that what is required is to approximate the true log likelihood uniformly, and that with increased accuracy of the approximation follows improved properties of the resulting approximate MLE. Uniform convergence on compacts holds in our case, because $l_n(\theta, x_0)$ is continuous in $\theta$, implied by the combination of so-called epi-convergence and hypoconvergence of an approximation $\tilde{l}_n(\theta, x_0)$ [see Geyer (1994), page 273]. Moreover, Geyer also proved that both of these modes of convergence can be obtained by an importance sampling approach, in which the unobserved states are simulated using MCMC under a fixed reference parameter [Geyer (1994), Theorem 2]. Of course, to obtain the required rate of convergence of the approximation, with increasing $n$ an increasing number of importance samples must be taken.

Approximation of the log likelihood using particle filters is described, for instance, in the above-mentioned paper by Pitt (2002), who also devised a method to smooth the approximation to a continuous function; this method works for univariate state variables only, however. At present we know of no formal proofs that particle filters approximate the true log likelihood uniformly, but strongly conjecture that they do under general assumptions.

8.3. *A numerical example.* We now turn to a specific numerical example, in which we shall employ an MCEM algorithm. Localization and tracking of narrow band moving sources by a passive array is one of the fundamental problems in radar, communication and sonar [see Ng, Larocque and Reilly (2001), Orton and Fitzgerald (2002) and references therein]. This problem can be stated as follows. Consider a uniform linear array of $d$ sensors receiving a narrowband signal from a far-field source with unknown time-varying direction of arrival (DOA). Under the classical narrowband array processing model the received signal at time $k$, the $d \times 1$ array observation vector $Y_k$, can be expressed as

$$W_k = W_{k-1} + \eta_k, \tag{27}$$

$$Y_k = S_k a(W_k) + \varepsilon_k, \tag{28}$$

where $a(w) = [1 e^{iw} \cdots e^{i(d-1)w}]^T$ is the $d \times 1$ *steering vector*, $S_k$ is the source waveform, $\eta_k$ is the state noise and $\varepsilon_k$ is the measurement noise. It is assumed that (i) $\{\eta_k\}$ are i.i.d. zero mean Gaussian with variance $\sigma_\eta^2$, (ii) $\{S_k\}$ are i.i.d. zero mean one-dimensional complex circular Gaussian, that is, $\mathbb{E} S_k = 0$ and $\mathbb{E} |S_k|^2 = \sigma_s^2$ and (iii) $\{\varepsilon_k\}$ are i.i.d. zero mean $d$-dimensional complex circular Gaussian, that is, $\mathbb{E} \varepsilon_k = 0$ and $\mathbb{E} \varepsilon_k \varepsilon_k^H = \sigma_\varepsilon^2 I_d$, where $x^H$ is the conjugate transpose (or Hermite transpose) of $x$ and $I_d$ is the $d \times d$ identity matrix. This is a hidden Markov model, or state space model, as there is no autoregression in the $Y$'s. We wish to estimate the parameter $\theta \triangleq (\sigma_\eta^2, \sigma_s^2, \sigma_\varepsilon^2)$ from the observed data $Y_1, \ldots, Y_n$.

Conditionally on the hidden variable $W_k$, $Y_k$ is a Gaussian complex vector with density $g_\theta(y_k|W_k)$, where

$$g_\theta(y|w) = \frac{1}{\pi^d \det \Sigma(w)} \exp\{-y^H \Sigma^{-1}(w) y\}$$

with

$$\Sigma(w) = \mathbb{E}[Y_k Y_k^H | W_k = w] = \sigma_s^2 a(w) a(w)^H + \sigma_\varepsilon^2 I_d.$$

It is easily checked that

$$\Sigma^{-1}(w) = -\frac{\sigma_s^2}{\sigma_\varepsilon^2 (d\sigma_s^2 + \sigma_\varepsilon^2)} a(w) a(w)^H + \frac{1}{\sigma_\varepsilon^2} I_d$$

and

$$\log g_\theta(y|w) = -d \log \pi - \log(\sigma_\varepsilon^{2(d-1)}(d\sigma_s^2 + \sigma_\varepsilon^2))$$
$$- \frac{1}{\sigma_\varepsilon^2} y^H y + \frac{\sigma_s^2}{\sigma_\varepsilon^2 (d\sigma_s^2 + \sigma_\varepsilon^2)} |a(w)^H y|^2.$$

Furthermore, with $r_\theta$ denoting the transition density of $\{W_k\}$,

$$\log r_\theta(w, w') = \log r_{\sigma_\eta^2}(w, w') = -\frac{1}{2} \log(2\pi \sigma_\eta^2) - \frac{1}{2\sigma_\eta^2}(w' - w)^2.$$

The above model is equivalent to an HMM on a compact state space. Indeed, identify the interval $[0, 2\pi)$ with the unit circle, which is a compact set, and put $X_k = W_k \bmod 2\pi$. It is then clear that $\{X_k\}$ is a Markov chain on $[0, 2\pi)$, with transition density $q_{\sigma_\eta^2}(x, x') = \sum_{\ell=-\infty}^{\infty} r_{\sigma_\eta^2}(x, x' + 2\pi\ell)$. The output density stays the same, that is, the conditional density of $Y_k$ given $X_k = x$ is $g_\theta(y|x)$. It is easily verified that the HMM $\{(X_k, Y_k)\}$ satisfies the regularity conditions in the previous sections.

Let $\theta$ and $\theta'$ denote two (potentially) different parameter values. The EM algorithm involves iterative maximization of the function $Q(\theta, \theta') = \mathbb{E}_\theta[\log p_{\theta'}(\mathbf{X}_1^n, \mathbf{Y}_1^n | X_0 = x_0) | \mathbf{Y}_1^n, X_0 = x_0]$. Specifically, if $\theta_p$ is the result of the $p$th iteration, then $\theta_{p+1}$ is the maximizer (in $\theta'$) of $Q(\theta_p, \theta')$, that is, $\theta_{p+1} = \arg\max_{\theta'} Q(\theta_p, \theta')$. For the present model, put $\beta(\theta) = \sum_{k=1}^n \mathbb{E}_\theta[|a(X_k)^H Y_k|^2 | \mathbf{Y}_1^n, X_0 = x_0]$. It is then straightforward to verify that the maximizer of the M-step of the EM algorithm is the triple $(\hat\sigma_\eta^2, \hat\sigma_s^2, \hat\sigma_\varepsilon^2)$ given by

(29) $$\hat\sigma_\eta^2 = \arg\max_v \mathbb{E}_{\theta_p}\left[\sum_{k=1}^n \log q_v(X_{k-1}, X_k) \bigg| \mathbf{Y}_1^n, X_0 = x_0\right],$$

$$\hat{\sigma}_s^2 = \frac{\beta(\theta_p) - \sum_{k=1}^n |Y_k|^2}{nd(d-1)}, \tag{30}$$

$$\hat{\sigma}_\varepsilon^2 = \frac{\sum_{k=1}^n |Y_k|^2 - \beta(\theta_p)/d}{n(d-1)}. \tag{31}$$

The conditional expectation $\beta(\theta)$ cannot be explicitly computed, let alone the expectation required to compute $\hat{\sigma}_\eta^2$. We note that we could also employ the representation with $\{W_k\}$ to simplify the implementation of this part of the M-step and the MCMC algorithm below, as there is then a sufficient statistic for the re-estimation of $\sigma_\eta^2$ as well, but this approach gives us less satisfying numerical results. We also note that although $q_\theta(x, x')$ is not available in closed form, it is straightforward to approximate it by a truncated sum as $r_{\sigma_\eta^2}(x, x' + 2\pi\ell)$ decays rapidly as $|\ell| \to \infty$.

In the MCEM approach, the conditional expectations above are replaced by sample means over a number of realizations of $\mathbf{X}_1^n$, conditional on $\mathbf{Y}_1^n$ and $X_0 = x_0$, obtained by Monte Carlo simulation. At each iteration $p$ we draw a sample of size $m_p$ of an $\mathbb{R}^n$-valued Markov chain $\{\hat{X}^{(\ell)}\}_{\ell \geq 0}$ with stationary distribution $\mathbb{P}_{\theta_p}(\mathbf{X}_1^n \in \cdot | \mathbf{Y}_1^n, X_0 = x_0)$. Many different solutions are available at this stage; in the simulations below, we use a random scan Metropolis–Hasting algorithm with transition kernel from $\hat{X}^{(\ell-1)} = \hat{x}$ to $\hat{X}^{(\ell)} = \hat{x}'$ defined in the following way:

1. Choose a time index $i$ uniformly on $\{1, \ldots, n\}$.
2. Simulate $\hat{x}_i'' \sim q_{\theta_p}(\hat{x}_{i-1}, \cdot)$.
3. Set $\hat{x}' = \hat{x}$ (these are $\mathbb{R}^n$-valued) and update the $i$th component of $\hat{x}'$, that is, $\hat{x}_i'$, to $\hat{x}_i''$ with probability

$$1 \wedge \frac{q_{\theta_p}(\hat{x}_{i-1}, \hat{x}_i'') q_{\theta_p}(\hat{x}_i'', \hat{x}_{i+1}) g_{\theta_p}(Y_i | \hat{x}_i'')}{q_{\theta_p}(\hat{x}_{i-1}, \hat{x}_i) q_{\theta_p}(\hat{x}_i, \hat{x}_{i+1}) g_{\theta_p}(Y_i | \hat{x}_i)} \times \frac{q_{\theta_p}(\hat{x}_{i-1}, \hat{x}_i)}{q_{\theta_p}(\hat{x}_{i-1}, \hat{x}_i'')}$$
$$= 1 \wedge \frac{q_{\theta_p}(\hat{x}_i'', \hat{x}_{i+1}) g_{\theta_p}(Y_i | \hat{x}_i'')}{q_{\theta_p}(\hat{x}_i, \hat{x}_{i+1}) g_{\theta_p}(Y_i | \hat{x}_i)}.$$

If $i = n$, this acceptance probability is modified to

$$1 \wedge \frac{q_{\theta_p}(\hat{x}_{i-1}, \hat{x}_i'') g_{\theta_p}(Y_i | \hat{x}_i'')}{q_{\theta_p}(\hat{x}_{i-1}, \hat{x}_i) g_{\theta_p}(Y_i | \hat{x}_i)} \times \frac{q_{\theta_p}(\hat{x}_{i-1}, \hat{x}_i)}{q_{\theta_p}(\hat{x}_{i-1}, \hat{x}_i'')} = 1 \wedge \frac{g_{\theta_p}(Y_i | \hat{x}_i'')}{g_{\theta_p}(Y_i | \hat{x}_i)}.$$

To guarantee convergence of the algorithm, the number of samples, $m_p$, should either be increased at each iteration or be selected in a data-driven manner at each iteration [see Booth and Hobert (1999) or Booth, Hobert and Jank (2001)]. For simplicity we did not implement such mechanisms but rather used a fixed large number of iterations at each step of the algorithm.

We simulated a single sample of size $n = 200$ from the model (27)–(28) with $d = 4$ and with the true value of the parameter $\theta = (\sigma_\eta^2, \sigma_s^2, \sigma_\varepsilon^2)$ being $\theta^* = (0.25, 0.64, 0.36)$. At each step of the MCEM procedure we generated a sample of size 40,000 by the random scan Metropolis–Hasting algorithm, after a burn-in of 20,000 iterations. The acceptance rate of the algorithm was about 40%. The re-estimation of $\sigma_\eta^2$ as in (29) was carried out by numerical optimization, and, in order to save computation time, of the total of 40,000 replications only every 400th was used for the corresponding sample average (i.e., 1,000 replications). The stationary distribution of $\{X_k\}$ is the uniform distribution on $[0, 2\pi)$, whence we fixed the initial state $x_0$ to its mean $\pi$. We remark that in this particular case the stationary distribution does not depend on $\theta$, whence it could have been employed in the algorithm. We started the MCEM algorithm from the true parameters as well as from four randomly chosen initial points for which each $\sigma^2$-parameter was drawn independently from a uniform distribution on $(0, 1)$. For each of the five initial points we ran the algorithm for 50 iterations. Figure 1 shows the trajectories for each initial point and parameter. Obviously, irrespective of the initial point the algorithm quickly finds the same approximation to the MLE, although the trajectories do not converge as the sample size $m_p$ in the algorithm stays bounded. The trajectories for $\sigma_\eta^2$ fluctuate a little more since, as described above, only 1,000 replications were used for its re-estimation.

Next we estimated the observed information, that is, the negative Hessian of the log likelihood. We departed from the missing information principle (22) and again replaced the expectations involved by sample means over simulated replications of $\mathbf{X}_1^n$ given $\mathbf{Y}_1^n$ and $X_0 = x_0$ obtained in the same way as above. Our approximation to the MLE, $\tilde{\theta}$ say, used for these computations was taken as the sample mean of the last 25 values of the trajectory obtained for the second randomly chosen starting point mentioned above; it was $\tilde{\theta} = (0.2793, 0.5756, 0.3466)$. After running the Metropolis–Hasting algorithm for a burn-in of 100,000 iterations we used another 200,000 iterations for the sample means. The resulting approximation of the observed information and its inverse were

$$\tilde{I} = \begin{pmatrix} 203.2 & -3.908 & 56.10 \\ -3.908 & 449.1 & 177.7 \\ 56.10 & 177.7 & 4169 \end{pmatrix},$$

$$\tilde{I}^{-1} = 10^{-3} \begin{pmatrix} 4.941 & 0.070 & -0.070 \\ 0.070 & 2.266 & -0.098 \\ -0.070 & -0.098 & 0.245 \end{pmatrix}.$$

The corresponding approximate 95% confidence intervals are $(0.1416, 0.4171)$, $(0.4823, 0.6689)$ and $(0.3159, 0.3773)$ for $\sigma_\eta^2$, $\sigma_s^2$ and $\sigma_\varepsilon^2$, respectively, and we see that they all cover the respective true values. We see that the variations in the MCEM estimates in Figure 1 are considerably smaller than the widths of the confidence intervals, which indicates that the MLE is well approximated and hence that the inverse observed information matrix is a good estimate of the covariance matrix of the approximate MLE as well. Obviously the widest interval is that for $\sigma_\eta^2$, which is not surprising, as this parameter is associated with the hidden state alone and hence, loosely speaking, "less observable" than the other ones. A simultaneous test for $H_0:\theta=\theta^*$ can be carried out by computing the test statistic $\chi^2 = (\tilde\theta-\theta^*)^T \tilde I (\tilde\theta-\theta^*)$, which approximately has a $\chi^2$ distribution with 3 degrees of freedom under the null hypothesis. We found $\chi^2 = 3.065$ and the corresponding $p$-value is 0.38. The null hypothesis could thus not be rejected.

## APPENDIX

### A.1. Proofs of technical lemmas.

PROOF OF LEMMA 3. Assume $m' \geq m$. Note that

$$\bar p_\theta(Y_k|\overline{\mathbf{Y}}_{-m}^{k-1}, X_{-m}=x) - \bar p_\theta(Y_k|\overline{\mathbf{Y}}_{-m'}^{k-1}, X_{-m'}=x')$$

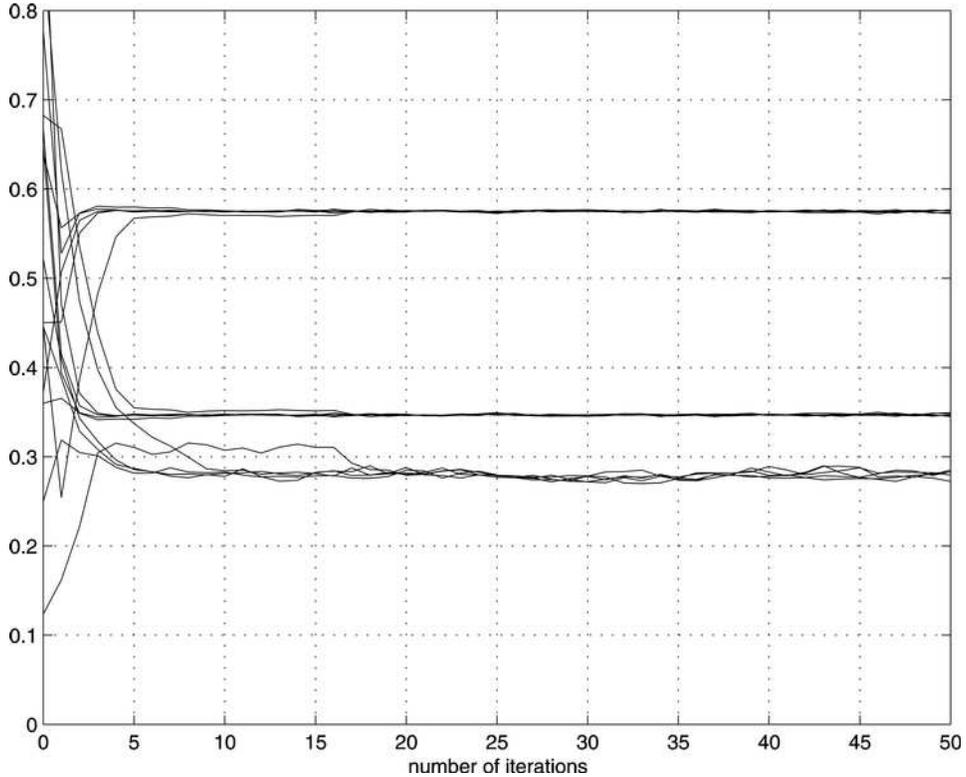

FIG. 1. *Convergence of the MCEM algorithm. Trajectories of the three parameters $\sigma_\eta^2$, $\sigma_s^2$ and $\sigma_\varepsilon^2$ for five runs of the MCEM algorithm, starting from five different initial points.*

$$= \iiint g_\theta(Y_k|\overline{\mathbf{Y}}_{k-1},x_k) q_\theta(x_{k-1},x_k) \mu(dx_k)$$

$$\times \overline{\mathbb{P}}_\theta(dx_{k-1}|X_{-m}=x_{-m}, \overline{\mathbf{Y}}_{-m}^{k-1}) \delta_x(dx_{-m})$$

$$- \iiint g_\theta(Y_k|\overline{\mathbf{Y}}_{k-1},x_k) q_\theta(x_{k-1},x_k) \mu(dx_k)$$

$$\times \overline{\mathbb{P}}_\theta(dx_{k-1}|X_{-m}=x_{-m}, \overline{\mathbf{Y}}_{-m}^{k-1}) \overline{\mathbb{P}}_\theta(dx_{-m}|\overline{\mathbf{Y}}_{-m'}^{k-1}, X_{-m'}=x').$$

Hence, by Corollary 1,

$$(32) \quad \begin{aligned} |\bar{p}_\theta(Y_k|\overline{\mathbf{Y}}_{-m}^{k-1}, X_{-m}=x) - \bar{p}_\theta(Y_k|\overline{\mathbf{Y}}_{-m'}^{k-1}, X_{-m'}=x')| \\ \leq \rho^{k+m-1} \sigma_+ \int g_\theta(Y_k|\overline{\mathbf{Y}}_{k-1},x) \mu(dx). \end{aligned}$$

Similarly we have

$$(33) \quad \begin{aligned} \bar{p}(Y_k|\overline{\mathbf{Y}}_{-m}^{k-1}, X_{-m}=x) \\ = \iint g_\theta(Y_k|\overline{\mathbf{Y}}_{k-1},x_k) q_\theta(x_{k-1},x_k) \mu(dx_k) \overline{\mathbb{P}}_\theta(dx_{k-1}|\overline{\mathbf{Y}}_{-m}^{k-1}, X_{-m}=x) \\ \geq \sigma_- \int g_\theta(Y_k|\overline{\mathbf{Y}}_{k-1},x) \mu(dx). \end{aligned}$$

The proof of (12) is concluded as in Lemma 2, and (13) follows by setting $m'=m$ and integrating w.r.t. $\overline{\mathbb{P}}_\theta(dx_{-m}|\overline{\mathbf{Y}}_{-m}^{k-1})$ in (32) and (33). To prove (14), notice that, by (33),

$$\sigma_- b_-(Y_k, \overline{\mathbf{Y}}_{k-1}) \leq \bar{p}_\theta(Y_k|\overline{\mathbf{Y}}_{-m}^{k-1}, X_{-m}=x) \leq b_+. \qquad \square$$

PROOF OF LEMMA 4. We will first prove that for any fixed $x \in \mathcal{X}$ and any $m$, $\Delta_{0,m,x}(\theta)$ is continuous w.r.t. $\theta$. We have

$$\bar{p}_\theta(Y_0|\overline{\mathbf{Y}}_{-m}^{-1}, X_{-m}=x) = \frac{\bar{p}_\theta(\mathbf{Y}_{-m+1}^0|\overline{\mathbf{Y}}_{-m}, X_{-m}=x)}{\bar{p}_\theta(\mathbf{Y}_{-m+1}^{-1}|\overline{\mathbf{Y}}_{-m}, X_{-m}=x)}$$

where, for $j \in \{-1, 0\}$,

$$(34) \quad \begin{aligned} \bar{p}_\theta(\mathbf{Y}_{-m+1}^j|\overline{\mathbf{Y}}_{-m}, X_{-m}=x) \\ = \int q_\theta(x, x_{-m+1}) \prod_{i=-m+2}^{j} q_\theta(x_{i-1}, x_i) \\ \times \prod_{i=-m+1}^{j} g_\theta(Y_i|\overline{\mathbf{Y}}_{i-1}, x_i) \mu^{\otimes(m+j)}(d\mathbf{x}_{-m+1}^j). \end{aligned}$$

Thus $\bar{p}_\theta(\mathbf{Y}^j_{-m+1}|\overline{\mathbf{Y}}_{-m}, X_{-m} = x)$ is continuous w.r.t. $\theta$ by continuity of $q_\theta$ and $g_\theta$ and the bounded convergence theorem; the integrand is bounded by $(\sigma_+ b_+)^{m+j}$. Since $\{\Delta_{0,m,x}(\theta)\}$ converges uniformly w.r.t. $\theta \in \Theta$, $\overline{\mathbb{P}}_{\theta^*}$-a.s., $\Delta_{0,\infty}(\theta)$ is continuous w.r.t. $\theta \in \Theta$, $\overline{\mathbb{P}}_{\theta^*}$-a.s., and the proof follows using Lemma 3 and the dominated convergence theorem. $\square$

PROOF OF PROPOSITION 2. By Lemma 2 it is sufficient to prove that
$$\limsup_{n\to\infty} \sup_{\theta\in\Theta} |n^{-1} l_n(\theta) - l(\theta)| = 0, \qquad \overline{\mathbb{P}}_{\theta^*}\text{-a.s.}$$

Furthermore, since $\Theta$ is compact, we only need to prove that for all $\theta \in \Theta$,
$$\limsup_{\delta\to 0} \limsup_{n\to\infty} \sup_{|\theta'-\theta|\leq\delta} |n^{-1} l_n(\theta') - l(\theta)| = 0, \qquad \overline{\mathbb{P}}_{\theta^*}\text{-a.s.}$$

Decompose the difference as
$$\limsup_{\delta\to 0} \limsup_{n\to\infty} \sup_{|\theta'-\theta|\leq\delta} |n^{-1} l_n(\theta') - l(\theta)|$$
$$= \limsup_{\delta\to 0} \limsup_{n\to\infty} \sup_{|\theta'-\theta|\leq\delta} |n^{-1} l_n(\theta') - n^{-1} l_n(\theta)|$$
$$\leq A + B + C,$$

where
$$A = \limsup_{\delta\to 0} \limsup_{n\to\infty} \sup_{|\theta'-\theta|\leq\delta} n^{-1} \sum_{k=1}^n |\Delta_{k,0}(\theta') - \Delta_{k,\infty}(\theta')|,$$
$$B = \limsup_{\delta\to 0} \limsup_{n\to\infty} \sup_{|\theta'-\theta|\leq\delta} n^{-1} \sum_{k=1}^n |\Delta_{k,\infty}(\theta') - \Delta_{k,\infty}(\theta)|,$$
$$C = \limsup_{n\to\infty} n^{-1} \sum_{k=1}^n |\Delta_{k,\infty}(\theta) - \Delta_{k,0}(\theta)|.$$

The terms $A$ and $C$ are zero by Corollary 2, and by the ergodic theorem and Lemma 4,
$$B \leq \limsup_{\delta\to 0} \limsup_{n\to\infty} n^{-1} \sum_{k=1}^n \sup_{|\theta'-\theta|\leq\delta} |\Delta_{k,\infty}(\theta') - \Delta_{k,\infty}(\theta)|$$
$$= \limsup_{\delta\to 0} \overline{\mathbb{E}}_{\theta^*}\left[\sup_{|\theta'-\theta|\leq\delta} |\Delta_{0,\infty}(\theta') - \Delta_{0,\infty}(\theta)|\right]$$
$$= 0, \qquad \overline{\mathbb{P}}_{\theta^*}\text{-a.s.} \qquad \square$$

PROOF OF LEMMA 5. We will show that for all $\ell > 0$,
$$(35) \qquad \sup_{i\leq 0} |\bar{p}_\theta(\mathbf{Y}^{k+\ell}_k | \overline{\mathbf{Y}}^0_{-i}) - \bar{p}_\theta(\mathbf{Y}^{k+\ell}_k)| \to 0, \qquad \overline{\mathbb{P}}_{\theta^*}\text{-a.s. as } k \to \infty.$$

By stationarity, this implies the statement of the lemma.

First recall that $z_s = (x_s, y_s, \ldots, y_1)$ and note that

$$\bar{p}_\theta(z_s|\bar{\mathbf{y}}^0_{-i}) = \iint \prod_{j=1}^s q_\theta(x_{j-1}, x_j) g_\theta(y_j|x_j, \bar{\mathbf{y}}_{j-1}) \overline{\mathbb{P}}_\theta(dx_0|\bar{\mathbf{y}}^0_{-i}) \mu^{\otimes(s-1)}(dx_1^{s-1})$$

$$\leq \sigma_+ \iint \prod_{j=2}^s q_\theta(x_{j-1}, x_j) \prod_{j=1}^s g_\theta(y_j|x_j, \bar{\mathbf{y}}_{j-1})$$

$$\times \overline{\mathbb{P}}_\theta(dx_0|\bar{\mathbf{y}}^0_{-i}) \mu^{\otimes(s-1)}(dx_1^{s-1})$$

$$= \sigma_+ h_\theta(z_s),$$

say, where $h_\theta(z_s)$ implicitly depends on $\bar{\mathbf{y}}_0$, but not on $i$, and integrates to unity (it is a density w.r.t. $\mu \otimes \bar{\nu}$). Furthermore,

$$|\bar{p}_\theta(\mathbf{Y}^{k+\ell}_k|\overline{\mathbf{Y}}^0_{-i}) - \bar{p}_\theta(\mathbf{Y}^{k+\ell}_k)|$$

$$\leq \iint \bar{p}_\theta(\mathbf{Y}^{k+\ell}_k|z_{k-1}) |\Pi^{k-s-1}(z_s, dz_{k-1}) - \pi_\theta(dz_{k-1})|$$

$$\times \bar{p}_\theta(z_s|\overline{\mathbf{Y}}^0_{-i})(\mu \otimes \bar{\nu})(dz_s)$$

$$\leq b_+^\ell \sigma_+ \int \|\Pi^{k-s-1}(z_s, \cdot) - \pi_\theta\|_{\mathrm{TV}} h_\theta(z_s) (\mu \otimes \bar{\nu})(dz_s);$$

the bound on $\bar{p}_\theta(\mathbf{Y}^{k+\ell}_k|z_{k-1})$ follows as in (34). Now (35) is a result of the above, (3) and dominated convergence. □

Let, for $0 \leq k \leq m$,

$$U_{k,m}(\theta) \triangleq \log \bar{p}_\theta(\mathbf{Y}^p_1|\overline{\mathbf{Y}}_0, \overline{\mathbf{Y}}^{-k}_{-m}), \qquad U(\theta) \triangleq \log \bar{p}_\theta(\mathbf{Y}^p_1|\overline{\mathbf{Y}}_0).$$

PROOF OF LEMMA 7. It is enough to show that, for all $\theta \in \Theta$,

(36) $$\lim_{k \to \infty} \overline{\mathbb{E}}_{\theta^*}\left[\sup_{m \geq k} |U_{k,m}(\theta) - U(\theta)|\right] = 0.$$

Put

$$A_{k,m} = \bar{p}_\theta(\mathbf{Y}^p_{-s+1}|\overline{\mathbf{Y}}^{-k}_{-m}), \qquad A = \bar{p}_\theta(\mathbf{Y}^p_{-s+1}),$$

$$B_{k,m} = \bar{p}_\theta(\overline{\mathbf{Y}}_0|\overline{\mathbf{Y}}^{-k}_{-m}), \qquad B = \bar{p}_\theta(\overline{\mathbf{Y}}_0).$$

Then

(37) $$|\bar{p}_\theta(\mathbf{Y}^p_1|\overline{\mathbf{Y}}_0, \overline{\mathbf{Y}}^{-k}_{-m}) - \bar{p}_\theta(\mathbf{Y}^p_1|\overline{\mathbf{Y}}_0)| = \left|\frac{A_{k,m}}{B_{k,m}} - \frac{A}{B}\right|$$

$$\leq \frac{B|A_{k,m} - A| + A|B_{k,m} - B|}{BB_{k,m}}.$$

By conditioning on $(X_{-s}, \overline{\mathbf{Y}}_{-s})$ [cf. (34)] and utilizing (A1)(b), it follows that

$$B \geq \sigma_-^s \int \prod_{i=-s+1}^{0} \int g_\theta(Y_i|Y_{i-1},\ldots,Y_{-s+1},y_{-s},\ldots,y_{i-s},x)$$
(38)
$$\times \mu(dx)\overline{\mathbb{P}}_\theta(\overline{\mathbf{Y}}_{-s} \in d\bar{\mathbf{y}}_{-s}) > 0, \qquad \overline{\mathbb{P}}_{\theta^*}\text{-a.s.}$$

Hence, by Lemma 5, with $\overline{\mathbb{P}}_{\theta^*}$-probability arbitrarily close to 1, $B_{k,m}(\omega)$ is uniformly bounded away from zero for $m \geq k$ and $k$ sufficiently large, and Lemma 5 and (37) show that

$$\lim_{k\to\infty} \sup_{m\geq k} |\bar{p}_\theta(\mathbf{Y}_1^p|\overline{\mathbf{Y}}_0,\overline{\mathbf{Y}}_{-m}^{-k}) - \bar{p}_\theta(\mathbf{Y}_1^p|\overline{\mathbf{Y}}_0)| = 0 \qquad \text{in } \overline{\mathbb{P}}_{\theta^*}\text{-probability.}$$

Using the inequality $|\log x - \log y| \leq |x-y|/(x\wedge y)$ and (38) once again, we find that

$$\lim_{k\to\infty} \sup_{m\geq k} |U_{k,m}(\theta) - U(\theta)| = 0 \qquad \text{in } \overline{\mathbb{P}}_{\theta^*}\text{-probability,}$$

and (36) follows using dominated convergence provided

$$\overline{\mathbb{E}}_{\theta^*}\left[\sup_k \sup_{m\geq k} |U_{k,m}(\theta)|\right] < \infty.$$

This expectation is indeed finite since $\bar{p}_\theta(\mathbf{Y}_1^p|\overline{\mathbf{Y}}_0,\overline{\mathbf{Y}}_{-m}^{-k})$ is bounded from below by $\sigma_-^p \prod_1^p b_-(\overline{\mathbf{Y}}_{i-1},Y_i)$ and from above by $(\sigma_+ b_+)^p$ [cf. (34)], and the logarithms of these bounds are in $L^1(\overline{\mathbb{P}}_{\theta^*})$. □

**A.2. Proof of Proposition 4.** We preface the proof with several lemmas. For convenience, Proposition 4 will be proved for $q=1$. Adaptations to general $q$ are obvious.

Define for $k \geq 1$, $m \geq 0$ and $x \in \mathcal{X}$,

$$\Delta_{k,m,x}(\theta) \triangleq \overline{\mathbb{E}}_\theta\left[\sum_{i=-m+1}^{k} \varphi(\theta,Z_i)\Big|\overline{\mathbf{Y}}_{-m}^k, X_{-m} = x\right]$$
$$- \overline{\mathbb{E}}_\theta\left[\sum_{i=-m+1}^{k-1} \varphi(\theta,Z_i)\Big|\overline{\mathbf{Y}}_{-m}^{k-1}, X_{-m} = x\right].$$

Along the same lines as in Lemma 9, for $m,n \geq 0$ and $0 < k \leq n+m-1$,

$$\overline{\mathbb{P}}_\theta(X_{n-k} \in A|\mathbf{X}_{n-k+1}^n, \overline{\mathbf{Y}}_{-m}^n, X_{-m} = x)$$
$$= \overline{\mathbb{P}}_\theta(X_{n-k} \in A|X_{n-k+1}, \overline{\mathbf{Y}}_{-m}^n, X_{-m} = x) \geq \frac{\sigma_-}{\sigma_+}\check{\mu}_k(\overline{\mathbf{Y}}_{-m}^{n-k}, X_{-m} = x, A),$$

where $\breve{\mu}_k(\overline{\mathbf{Y}}_{-m}^{n-k}, X_{-m} = x, \cdot)$ is a probability measure. The result above in particular implies that

$$\|\overline{\mathbb{P}}_\theta(X_i \in \cdot | \overline{\mathbf{Y}}_{-m}^n, X_{-m} = x) - \overline{\mathbb{P}}_\theta(X_i \in \cdot | \overline{\mathbf{Y}}_{-m}^{n-1}, X_{-m} = x)\|_{\mathrm{TV}} \leq \rho^{n-i-1}.$$
(39)

LEMMA 13. *Under the assumptions of Proposition 4 there exists a random variable $K \in L^1(\overline{\mathbb{P}}_{\theta^*})$ such that, for all $k \geq 1$ and $0 \leq m \leq m'$,*

(40) $\sup_{x \in \mathcal{X}} \sup_{\theta \in G} |\Delta_{k,m,x}(\theta) - \Delta_{k,m}(\theta)| \leq K(k \vee m)^2 \rho^{(k+m)/2}, \qquad \overline{\mathbb{P}}_{\theta^*}\text{-}a.s.,$

(41) $\sup_{x \in \mathcal{X}} \sup_{\theta \in G} |\Delta_{k,m,x}(\theta) - \Delta_{k,m',x}(\theta)| \leq K(k \vee m)^2 \rho^{(k+m)/2}, \qquad \overline{\mathbb{P}}_{\theta^*}\text{-}a.s.$

PROOF. The proof is along the same lines as the proof of Lemma 10, using (39). Put $\|\varphi_i\|_\infty = \sup_{x \in \mathcal{X}} \sup_{\theta \in G} |\varphi(\theta, x, Y_i)|$. Combining the relations

$$|\overline{\mathbb{E}}_\theta[\varphi(\theta, Z_i) | \overline{\mathbf{Y}}_{-m}^k, X_{-m} = x] - \overline{\mathbb{E}}_\theta[\varphi(\theta, Z_i) | \overline{\mathbf{Y}}_{-m}^k]| \leq 2\|\varphi_i\|_\infty \rho^{i+m},$$

$$|\overline{\mathbb{E}}_\theta[\varphi(\theta, Z_i) | \overline{\mathbf{Y}}_{-m}^k, X_{-m} = x] - \overline{\mathbb{E}}_\theta[\varphi(\theta, Z_i) | \overline{\mathbf{Y}}_{-m}^{k-1}, X_{-m} = x]|$$
$$\leq 2\|\varphi_i\|_\infty \rho^{k-i-1},$$

$$|\overline{\mathbb{E}}_\theta[\varphi(\theta, Z_i) | \overline{\mathbf{Y}}_{-m}^k] - \overline{\mathbb{E}}_\theta[\varphi(\theta, Z_i) | \overline{\mathbf{Y}}_{-m}^{k-1}]| \leq 2\|\varphi_i\|_\infty \rho^{k-i-1},$$

we obtain

$$|\Delta_{k,m,x}(\theta) - \Delta_{k,m}(\theta)|$$

$$\leq 4 \sum_{i=-m+1}^k \|\varphi_i\|_\infty (\rho^{i+m} \wedge \rho^{k-1-i})$$

$$\leq 4 \max_{-m \leq i \leq k} \|\varphi_i\|_\infty \sum_{i=-m+1}^k (\rho^{i+m} \wedge \rho^{k-1-i})$$

$$\leq 4 \sum_{i=-m}^k (|i| \vee 1)^2 \frac{1}{(|i| \vee 1)^2} \|\varphi_i\|_\infty \left( \sum_{i \leq (k-m-1)/2} \rho^{k-1-i} + \sum_{i \geq (k-m-1)/2} \rho^{i+m} \right)$$

$$\leq 8(k \vee m)^2 \sum_{i=-\infty}^\infty \frac{1}{(|i| \vee 1)^2} \|\varphi_i\|_\infty \frac{\rho^{(k+m-1)/2}}{1-\rho},$$

which proves the first part of the lemma.

For the second part we also use the bound

$$|\overline{\mathbb{E}}_\theta[\varphi(\theta, Z_i) | \overline{\mathbf{Y}}_{-m}^k, X_{-m} = x] - \overline{\mathbb{E}}_\theta[\varphi(\theta, Z_i) | \overline{\mathbf{Y}}_{-m'}^k, X_{-m'} = x]|$$
$$\leq 2\|\varphi_i\|_\infty \rho^{i+m \wedge m'}$$

to obtain

$$|\Delta_{k,m,x}(\theta) - \Delta_{k,m',x}(\theta)|$$
$$\leq 4 \sum_{i=-m+1}^{k} \|\varphi_i\|_\infty (\rho^{i+m} \wedge \rho^{k-1-i}) + 2 \sum_{i=-m'+1}^{-m} \|\varphi_i\|_\infty \rho^{k-1-i}.$$

Here the first term on the right-hand side is bounded as above. Since $-i/2 \leq (k-m)/2 - i$ for $i \leq -m$, the second term can be bounded as

$$2 \sum_{i=-m'+1}^{-m} \|\varphi_i\|_\infty \rho^{k-1-i} \leq 2\rho^{(k+m)/2} \sum_{i=-m'+1}^{-m} \|\varphi_i\|_\infty \rho^{(k-m)/2-1-i}$$

$$\leq 2\rho^{(k+m)/2} \sum_{i=-m'+1}^{-m} \|\varphi_i\|_\infty \rho^{-i/2-1}$$

$$\leq 2\rho^{(k+m)/2} \sum_{i=-\infty}^{\infty} \|\varphi_i\|_\infty \rho^{|i|/2-1}$$

and the proof is complete. □

By Lemma 13, for all $x \in \mathcal{X}$ and $k \geq 1$, $\{\Delta_{k,m,x}(\theta)\}_{m \geq 0}$ converges uniformly w.r.t. $\theta \in G$ $\overline{\mathbb{P}}_{\theta^*}$-a.s. and in $L^1(\overline{\mathbb{P}}_{\theta^*})$ to a random variable that we denote by $\Delta_{k,\infty}(\theta)$; by (40) this limit does not depend on $x$. Lemma 13 also immediately implies that

$$n^{-1} \sum_{k=1}^{n} \sup_{\theta \in G} |\Delta_{k,0}(\theta) - \Delta_{k,\infty}(\theta)| \to 0, \qquad \overline{\mathbb{P}}_{\theta^*}\text{-a.s. and in } L^1(\overline{\mathbb{P}}_{\theta^*}).$$

LEMMA 14. *Under the assumptions of Proposition 4, for all $x \in \mathcal{X}$ and $m \geq 0$ the function $\theta \mapsto \Delta_{0,m,x}(\theta)$ is $\overline{\mathbb{P}}_{\theta^*}$-a.s. continuous on $G$. In addition, for all $\theta \in G$ and all $x \in \mathcal{X}$,*

$$\lim_{\delta \to 0} \overline{\mathbb{E}}_{\theta^*} \left[ \sup_{|\theta'-\theta| \leq \delta} |\Delta_{0,m,x}(\theta') - \Delta_{0,m,x}(\theta)| \right] = 0.$$

PROOF. Note that $|\Delta_{0,m,x}(\theta)| \leq 2\sum_{i=-m+1}^{0} \|\varphi_i\|_\infty$. Thus, under the assumptions of Proposition 4, $\Delta_{0,m,x}(\theta)$ is uniformly bounded w.r.t. $\theta$ by a random variable in $L^1(\overline{\mathbb{P}}_{\theta^*})$. It hence suffices to show that for $-m < i \leq 0$,

$$\lim_{\delta \to 0} \sup_{|\theta'-\theta| \leq \delta} |\overline{\mathbb{E}}_{\theta'}[\varphi(\theta', Z_i)|\overline{\mathbf{Y}}_{-m}^0, X_{-m} = x]$$
$$- \overline{\mathbb{E}}_{\theta}[\varphi(\theta, Z_i)|\overline{\mathbf{Y}}_{-m}^0, X_{-m} = x]| = 0, \qquad \overline{\mathbb{P}}_{\theta^*}\text{-a.s.}$$

Write

(42)
$$\overline{\mathbb{E}}_\theta[\varphi(\theta, Z_i)|\overline{\mathbf{Y}}^0_{-m}, X_{-m} = x]$$
$$= \int \varphi(\theta, x_i, Y_i)\bar{p}_\theta(X_i = x_i|\overline{\mathbf{Y}}^0_{-m}, X_{-m} = x)\mu(dx_i)$$

and note that for all $x_i$, $\varphi(\theta, x_i, Y_i)$ is continuous w.r.t. $\theta$ and that this factor is bounded by $\|\varphi_i\|_\infty < \infty$. Moreover,

$$\bar{p}_\theta(X_i = x_i|\overline{\mathbf{Y}}^0_{-m}, X_{-m} = x) = \frac{\bar{p}_\theta(X_i = x_i, \mathbf{Y}^0_{-m+1}|\overline{\mathbf{Y}}_{-m}, X_{-m} = x)}{\bar{p}_\theta(\mathbf{Y}^0_{-m+1}|\overline{\mathbf{Y}}_{-m}, X_{-m} = x)}.$$

Here $\bar{p}_\theta(\mathbf{Y}^0_{-m+1}|\overline{\mathbf{Y}}_{-m}, X_{-m} = x)$ is continuous w.r.t. $\theta$ (see the proof of Lemma 4), and using (34) we find that this density is bounded from below by

$$\sigma_-^m \prod_{i=-m+1}^{0} \int g_\theta(Y_i|\overline{\mathbf{Y}}_{i-1}, x_i)\mu(dx_i) > 0$$

uniformly w.r.t. $\theta$. In a similar fashion $\bar{p}_\theta(X_i = x_i, \mathbf{Y}^0_{-m+1}|\overline{\mathbf{Y}}_{-m}, X_{-m} = x)$ is continuous in $\theta$ and bounded from above by $(\sigma_+ b_+)^m$. We conclude that $\bar{p}_\theta(X_i = x_i|\overline{\mathbf{Y}}_{-m}, X_{-m} = x)$ is continuous in $\theta$ and bounded from above uniformly w.r.t. $\theta$. Hence the integrand in (42) is continuous in $\theta$ and bounded from above uniformly w.r.t. $\theta$. Dominated convergence shows that the left-hand side of (42) is continuous in $\theta$ and the proof is complete. □

By Lemma 13 $\Delta_{0,m,x}(\theta)$ is a uniform Cauchy sequence w.r.t. $\theta$ $\overline{\mathbb{P}}_{\theta^*}$-a.s. and in $L^1(\overline{\mathbb{P}}_{\theta^*})$, and by Lemma 14 $\Delta_{0,m,x}(\theta)$ is continuous w.r.t. $\theta$ on $G$ $\overline{\mathbb{P}}_{\theta^*}$-a.s. and in $L^1(\overline{\mathbb{P}}_{\theta^*})$ for each $m$. Hence it follows that $\Delta_{0,\infty}(\theta)$ is continuous w.r.t. $\theta$ on $G$ $\overline{\mathbb{P}}_{\theta^*}$-a.s. and in $L^1(\overline{\mathbb{P}}_{\theta^*})$, that is, for each $\theta \in G$,

(43) $\lim_{\delta \to 0} \sup_{|\theta'-\theta|\leq\delta} |\Delta_{0,\infty}(\theta') - \Delta_{0,\infty}(\theta)| = 0, \qquad \overline{\mathbb{P}}_{\theta^*}$-a.s. and in $L^1(\overline{\mathbb{P}}_{\theta^*})$.

REMARK 8. It is important to stress at this point that the result above *does not imply* that $\Delta_{0,m}(\theta)$ is continuous w.r.t. $\theta$ because, contrary to JP, we do not assume any kind of regularity condition for the stationary distribution as a function of $\theta$. Nevertheless, we have proved above that $\Delta_{0,\infty}(\theta)$ is continuous.

We may now prove a locally uniform law of large numbers.

LEMMA 15. *Under the assumptions of Proposition 4, for all $\theta \in G$,*

$$\lim_{\delta \to 0} \lim_{n \to \infty} \sup_{|\theta'-\theta|\leq\delta} \left| n^{-1} \sum_{k=1}^{n} \Delta_{k,\infty}(\theta') - \overline{\mathbb{E}}_{\theta^*}[\Delta_{0,\infty}(\theta)] \right| = 0, \qquad \overline{\mathbb{P}}_{\theta^*}\text{-a.s.}$$

PROOF. Write

$$\sup_{|\theta'-\theta|\leq\delta}\left|n^{-1}\sum_{k=1}^{n}\Delta_{k,\infty}(\theta') - \overline{\mathbb{E}}_{\theta^*}[\Delta_{0,\infty}(\theta)]\right|$$

$$\leq \sup_{|\theta'-\theta|\leq\delta}\left|n^{-1}\sum_{k=1}^{n}(\Delta_{k,\infty}(\theta') - \Delta_{k,\infty}(\theta))\right|$$

$$+ \left|n^{-1}\sum_{k=1}^{n}\Delta_{k,\infty}(\theta) - \overline{\mathbb{E}}_{\theta^*}[\Delta_{0,\infty}(\theta)]\right|$$

$$\leq n^{-1}\sum_{k=1}^{n}\sup_{|\theta'-\theta|\leq\delta}|\Delta_{k,\infty}(\theta') - \Delta_{k,\infty}(\theta)|$$

$$+ \left|n^{-1}\sum_{k=1}^{n}\Delta_{k,\infty}(\theta) - \overline{\mathbb{E}}_{\theta^*}[\Delta_{0,\infty}(\theta)]\right|.$$

As $n \to \infty$, the first term on the right-hand side tends to

$$\overline{\mathbb{E}}_{\theta^*}\left[\sup_{|\theta'-\theta|\leq\delta}|\Delta_{0,\infty}(\theta') - \Delta_{0,\infty}(\theta)|\right], \qquad \overline{\mathbb{P}}_{\theta^*}\text{-a.s.},$$

an expression which, by (43), vanishes when $\delta \to 0$. The second term vanishes $\overline{\mathbb{P}}_{\theta^*}$-a.s. as $n \to \infty$ by the ergodic theorem. This completes the proof. $\square$

We have now at hand all the necessary elements to prove Proposition 4.

PROOF OF PROPOSITION 4. Convergence of $\Delta_{k,m}(\theta)$ and continuity of $\overline{\mathbb{E}}_{\theta^*}[\Delta_{0,\infty}(\theta)]$ have been proved above, so it remains to show the last part of the proposition.

Note that

$$\mathbb{E}_\theta\left[\sum_{i=1}^{n}\varphi(\theta, \mathbf{Z}_{i-q+1}^i)\Big|\overline{\mathbf{Y}}_0^n, X_0 = x_0\right] = \sum_{k=1}^{n}\Delta_{k,0,x_0}(\theta).$$

Letting $m' \to \infty$ in Lemma 13 we find that $|\Delta_{k,0,x_0}(\theta) - \Delta_{k,\infty}(\theta)| \leq Kk^2\rho^{k/2}$ $\overline{\mathbb{P}}_{\theta^*}$-a.s. and hence it is sufficient to prove that

$$\lim_{\delta\to 0}\lim_{n\to\infty}\sup_{|\theta'-\theta|\leq\delta}\left|n^{-1}\sum_{k=1}^{n}\Delta_{k,\infty}(\theta') - \overline{\mathbb{E}}_{\theta^*}[\Delta_{0,\infty}(\theta)]\right| = 0, \qquad \overline{\mathbb{P}}_{\theta^*}\text{-a.s.}$$

This, however, is Lemma 15. $\square$

**A.3. Proof of Proposition 5.** The proof of Proposition 5 closely follows the proof of Proposition 4. Only the main adaptations from the proof are presented. We gather in the following lemma some of the required bounds for the conditional covariance. In the proof of Proposition 5 we will consider for convenience $q = 1$, and we let $\phi_{\theta,i} \triangleq \phi(\theta, Z_i)$ and $\|\phi_i\|_\infty \triangleq \sup_{\theta \in G} \sup_{x \in \mathcal{X}} |\phi(\theta, x, Y_i)|$.

LEMMA 16. *Under the assumptions of Proposition 5, for all $m' \geq m \geq 0$, all $-m < i, j \leq n$, all $\theta \in G$ and all $x \in \mathcal{X}$,*

$$|\overline{\text{cov}}_\theta[\phi_{\theta,i}, \phi_{\theta,j}|\overline{\mathbf{Y}}^n_{-m}]| \leq 2\rho^{|i-j|}\|\phi_i\|_\infty \|\phi_j\|_\infty,$$

$$|\overline{\text{cov}}_\theta[\phi_{\theta,i}, \phi_{\theta,j}|\overline{\mathbf{Y}}^n_{-m}, X_{-m} = x]| \leq 2\rho^{|i-j|}\|\phi_i\|_\infty \|\phi_j\|_\infty,$$

$$|\overline{\text{cov}}_\theta[\phi_{\theta,i}, \phi_{\theta,j}|\overline{\mathbf{Y}}^n_{-m}, X_{-m} = x] - \overline{\text{cov}}_\theta[\phi_{\theta,i}, \phi_{\theta,j}|\overline{\mathbf{Y}}^n_{-m}]|$$
$$\leq 6\|\phi_i\|_\infty \|\phi_j\|_\infty \rho^{m+i \wedge j},$$

$$|\overline{\text{cov}}_\theta[\phi_{\theta,i}, \phi_{\theta,j}|\overline{\mathbf{Y}}^n_{-m}] - \overline{\text{cov}}_\theta[\phi_{\theta,i}, \phi_{\theta,j}|\overline{\mathbf{Y}}^{n+1}_{-m}]| \leq 6\|\phi_i\|_\infty \|\phi_j\|_\infty \rho^{n-i \vee j},$$

$$|\overline{\text{cov}}_\theta[\phi_{\theta,i}, \phi_{\theta,j}|\overline{\mathbf{Y}}^n_{-m}, X_{-m} = x] - \overline{\text{cov}}_\theta[\phi_{\theta,i}, \phi_{\theta,j}|\overline{\mathbf{Y}}^{n+1}_{-m}, X_{-m} = x]|$$
$$\leq 6\|\phi_i\|_\infty \|\phi_j\|_\infty \rho^{n-i \vee j}.$$

All these relations stem from Corollary 1, Lemma 9, (39) and observations such as, for $i < j$,

$$|\overline{\mathbb{P}}_\theta(X_i \in A, X_j \in B|\overline{\mathbf{Y}}^n_{-m}, X_{-m} = x)$$
$$- \overline{\mathbb{P}}_\theta(X_i \in A|\overline{\mathbf{Y}}^n_{-m}, X_{-m} = x)\overline{\mathbb{P}}_\theta(X_j \in B|\overline{\mathbf{Y}}^n_{-m}, X_{-m} = x)|$$
$$= \overline{\mathbb{P}}_\theta(X_i \in A|\overline{\mathbf{Y}}^n_{-m}, X_{-m} = x)$$
$$\times |\overline{\mathbb{P}}_\theta(X_j \in B|\overline{\mathbf{Y}}^n_{-m}, X_i \in A, X_{-m} = x) - \overline{\mathbb{P}}_\theta(X_j \in B|\overline{\mathbf{Y}}^n_{-m}, X_{-m} = x)|$$
$$\leq \rho^{j-i}.$$

Details of the proof are omitted for brevity.

For $x \in \mathcal{X}$ define

$$\Gamma_{k,m,x}(\theta) \triangleq \overline{\text{var}}_\theta\left[\sum_{i=-m+1}^k \phi_{\theta,i}\bigg|\overline{\mathbf{Y}}^k_{-m}, X_{-m} = x\right]$$
$$- \overline{\text{var}}_\theta\left[\sum_{i=-m+1}^{k-1} \phi_{\theta,i}\bigg|\overline{\mathbf{Y}}^{k-1}_{-m}, X_{-m} = x\right].$$

We again follow the pattern of proof consisting of showing that for each $k$ and $x \in \mathcal{X}$, the sequence $\{\Gamma_{k,m,x}(\theta)\}_{m \geq 0}$ is a uniform (w.r.t. $\theta \in G$) Cauchy sequence that converges to a limit which does not depend on $x$.

LEMMA 17. *Under the assumptions of Proposition 5 there exists a random variable $K \in L^1(\overline{\mathbb{P}}_{\theta^*})$ such that, for all $k \geq 1$ and $0 \leq m \leq m'$,*

$$(44) \quad \sup_{x \in \mathcal{X}} \sup_{\theta \in G} |\Gamma_{k,m,x}(\theta) - \Gamma_{k,m}(\theta)| \leq K(m+k)^3 \rho^{(k+m)/4}, \qquad \overline{\mathbb{P}}_{\theta^*}\text{-}a.s.,$$

$$(45) \sup_{x \in \mathcal{X}} \sup_{\theta \in G} |\Gamma_{k,m,x}(\theta) - \Gamma_{k,m',x}(\theta)| \leq K(m+k)^3 \rho^{(k+m)/8}, \qquad \overline{\mathbb{P}}_{\theta^*}\text{-}a.s.$$

PROOF. Let, for $a \leq b$, $S_a^b \triangleq \sum_{i=a}^{b} \phi_{\theta,i}$ (the dependence on $\theta$ is implicit). The difference $\Gamma_{k,m,x}(\theta) - \Gamma_{k,m}(\theta)$ may be decomposed as $A + 2B + C$, where

$$A = \overline{\mathrm{var}}_\theta[S_{-m+1}^{k-1} | \overline{\mathbf{Y}}_{-m}^k, X_{-m} = x] - \overline{\mathrm{var}}_\theta[S_{-m+1}^{k-1} | \overline{\mathbf{Y}}_{-m}^{k-1}, X_{-m} = x]$$
$$\quad - \overline{\mathrm{var}}_\theta[S_{-m+1}^{k-1} | \overline{\mathbf{Y}}_{-m}^k] + \overline{\mathrm{var}}_\theta[S_{-m+1}^{k-1} | \overline{\mathbf{Y}}_{-m}^{k-1}],$$
$$B = \overline{\mathrm{cov}}_\theta[S_{-m+1}^{k-1}, \phi_{\theta,k} | \overline{\mathbf{Y}}_{-m}^k, X_{-m} = x] - \overline{\mathrm{cov}}_\theta[S_{-m+1}^{k-1}, \phi_{\theta,k} | \overline{\mathbf{Y}}_{-m}^k],$$
$$C = \overline{\mathrm{var}}_\theta[\phi_{\theta,k} | \overline{\mathbf{Y}}_{-m}^k, X_{-m} = x] - \overline{\mathrm{var}}_\theta[\phi_{\theta,k} | \overline{\mathbf{Y}}_{-m}^k].$$

By applying Lemma 16, it follows that

$$|A| \leq 2 \sum_{-m+1 \leq i \leq j \leq k-1} (2 \times 6\rho^{m+i} \wedge 4 \times 2\rho^{j-i} \wedge 2 \times 6\rho^{k-j-1})$$
$$\times \max_{-m+1 \leq i \leq j \leq k-1} \|\phi_i\|_\infty \|\phi_j\|_\infty.$$

The Cauchy–Schwarz inequality yields

$$(46) \quad \begin{aligned} \max_{-m+1 \leq i \leq j \leq k} \|\phi_i\|_\infty \|\phi_j\|_\infty &\leq \left( \sum_{i=-m}^{k} \|\phi_i\|_\infty \right)^2 \\ &\leq \sum_{i=-m}^{k} (|i| \vee 1)^2 \sum_{i=-m}^{k} \frac{1}{(|i| \vee 1)^2} \|\phi_i\|_\infty^2 \\ &\leq (m^3 + k^3) \sum_{i=-\infty}^{\infty} \frac{1}{(|i| \vee 1)^2} \|\phi_i\|_\infty^2, \end{aligned}$$

where the last sum is in $L^1(\overline{\mathbb{P}}_{\theta^*})$. Furthermore, for $n \geq 0$,

$$\sum_{0 \leq i \leq j \leq n} (\rho^i \wedge \rho^{j-i} \wedge \rho^{n-j}) \leq 2 \sum_{0 \leq i \leq n/2} \sum_{i \leq j \leq n-i} (\rho^{n-j} \wedge \rho^{j-i})$$
$$\leq 2 \sum_{0 \leq i \leq n/2} \left( \sum_{i \leq j \leq (n+i)/2} \rho^{n-j} + \sum_{(n+i)/2 \leq j \leq n-i} \rho^{j-i} \right)$$
$$\leq \frac{4}{1-\rho} \sum_{0 \leq i \leq n/2} \rho^{(n-i)/2}$$

$$\leq \frac{4\rho^{n/4}}{(1-\rho)(1-\rho^{1/2})}.$$

This shows that $|A|$ is bounded by an expression as in the first part of the lemma.

Similarly we have

$$|B| \leq 6 \sum_{i=-m+1}^{k-1} (\rho^{m+i} \wedge \rho^{k-i}) \max_{-m+1 \leq i \leq k-1} \|\phi_i\|_\infty \|\phi_k\|_\infty.$$

For the maximum we can use the bound (46), and for the sum we note that, for $n \geq 0$,

$$\sum_{i=0}^{n} (\rho^i \wedge \rho^{n-i}) = \sum_{0 \leq i \leq n/2} \rho^{n-i} + \sum_{n/2 \leq i \leq n} \rho^i \leq \frac{2\rho^{n/2}}{1-\rho}.$$

Thus $|B|$ is bounded by an expression as in the first part of the lemma.

For $C$ we have $|C| \leq 6\rho^{k+m}\|\phi_k\|_\infty^2$, and the proof of the first part of the lemma is complete.

The difference $\Gamma_{k,m,x}(\theta) - \Gamma_{k,m',x}(\theta)$ may be decomposed as $A + 2B + C + D + 2E + 2F$, where

$$A = \overline{\mathrm{var}}_\theta[S_{-m+1}^{k-1}|\overline{\mathbf{Y}}_{-m}^{k}, X_{-m}=x] - \overline{\mathrm{var}}_\theta[S_{-m+1}^{k-1}|\overline{\mathbf{Y}}_{-m}^{k-1}, X_{-m}=x]$$
$$\quad - \overline{\mathrm{var}}_\theta[S_{-m+1}^{k-1}|\overline{\mathbf{Y}}_{-m'}^{k}, X_{-m'}=x] + \overline{\mathrm{var}}_\theta[S_{-m+1}^{k-1}|\overline{\mathbf{Y}}_{-m'}^{k-1}, X_{-m'}=x],$$

$$B = \overline{\mathrm{cov}}_\theta[S_{-m+1}^{k-1}, \phi_{\theta,k}|\overline{\mathbf{Y}}_{-m}^{k}, X_{-m}=x] - \overline{\mathrm{cov}}_\theta[S_{-m+1}^{k-1}, \phi_{\theta,k}|\overline{\mathbf{Y}}_{-m'}^{k}, X_{-m'}=x],$$

$$C = \overline{\mathrm{var}}_\theta[\phi_{\theta,k}|\overline{\mathbf{Y}}_{-m}^{k}, X_{-m}=x] - \overline{\mathrm{var}}_\theta[\phi_{\theta,k}|\overline{\mathbf{Y}}_{-m'}^{k}, X_{-m'}=x],$$

$$D = \overline{\mathrm{var}}_\theta[S_{-m'+1}^{-m}|\overline{\mathbf{Y}}_{-m'}^{k}, X_{-m'}=x] - \overline{\mathrm{var}}_\theta[S_{-m'+1}^{-m}|\overline{\mathbf{Y}}_{-m'}^{k-1}, X_{-m'}=x],$$

$$E = \overline{\mathrm{cov}}_\theta[S_{-m+1}^{k-1}, S_{-m'+1}^{-m}|\overline{\mathbf{Y}}_{-m'}^{k}, X_{-m'}=x]$$
$$\quad - \overline{\mathrm{cov}}_\theta[S_{-m+1}^{k-1}, S_{-m'+1}^{-m}|\overline{\mathbf{Y}}_{-m'}^{k-1}, X_{-m'}=x],$$

$$F = \overline{\mathrm{cov}}_\theta[S_{-m'+1}^{-m}, \phi_{\theta,k}|\overline{\mathbf{Y}}_{-m'}^{k}, X_{-m'}=x].$$

Here $|A|$, $|B|$ and $|C|$ can be bounded as above, using variants of the bounds in Lemma 16.

Before proceeding, we note that for $k \geq 1$, $m \geq 0$ and $i \leq 0$, the following implications hold:

$$\text{if } j \leq (k+i-1)/2, \quad \text{then } (|j|-1)/2 \leq (3k+i-3)/4 - j,$$
$$\text{if } (k+i-1)/2 \leq j \leq k-1, \quad \text{then } (|j|-1)/4 \leq j + (-k-3i+1)/4,$$

MLE FOR SWITCHING AUTOREGRESSION 49$$\text{if } i \leq -m, \qquad \text{then } |i|/8 \leq (k-2i-m)/8,$$
$$\text{if } i \leq -m, \qquad \text{then } 3|i|/4 \leq (k-m)/4 - i.$$

Using these inequalities, we can bound $|D|$ as

$$|D| \leq 2 \sum_{-m'+1 \leq i \leq j \leq -m} (6\rho^{k-1-j} \wedge 2 \times 2\rho^{j-i}) \|\phi_i\|_\infty \|\phi_j\|_\infty$$

$$\leq 12\rho^{(k+m-2)/8} \sum_{i=-m'+1}^{-m} \rho^{(k-2i-m)/8} \|\phi_i\|_\infty$$

$$\times \left( \sum_{i \leq j \leq (k+i-1)/2} \rho^{(3k+i-3)/4-j} \|\phi_j\|_\infty \right.$$

$$\left. + \sum_{(k+i-1)/2 < j \leq -m} \rho^{j+(-k-3i+1)/4} \|\phi_j\|_\infty \right)$$

$$\leq 12\rho^{(k+m-2)/8} \sum_{i=-\infty}^{\infty} \rho^{|i|/8} \|\phi_i\|_\infty \sum_{j=-\infty}^{\infty} \rho^{(|j|-1)/4} \|\phi_j\|_\infty.$$

By the assumptions, the right-hand side has the required form.

Similarly,

$$|E| \leq \sum_{i=-m'+1}^{-m} \sum_{j=-m+1}^{k-1} (6\rho^{k-1-j} \wedge 2 \times 2\rho^{j-i}) \|\phi_i\|_\infty \|\phi_j\|_\infty$$

$$\leq 6\rho^{(k+m-2)/8} \sum_{i=-m'+1}^{-m} \rho^{(k-2i-m)/8} \|\phi_i\|_\infty$$

$$\times \left( \sum_{-m+1 \leq j \leq (k+i-1)/2} \rho^{(3k+i-3)/4-j} \|\phi_j\|_\infty \right.$$

$$\left. + \sum_{(k+i-1)/2 < j \leq k-1} \rho^{j+(-k-3i+1)/4} \|\phi_j\|_\infty \right)$$

$$\leq 6\rho^{(k+m-2)/8} \sum_{i=-\infty}^{\infty} \rho^{|i|/8} \|\phi_i\|_\infty \sum_{j=-\infty}^{\infty} \rho^{(|j|-1)/4} \|\phi_j\|_\infty$$

and

$$|F| \leq \sum_{i=-m'+1}^{-m} 2\rho^{k-i} \|\phi_i\|_\infty \|\phi_k\|_\infty$$

$$= 2\rho^{(k+m)/4} \sum_{i=-m'+1}^{-m} \rho^{(k-m)/4-i} \|\phi_i\|_\infty \rho^{k/2} \|\phi_k\|_\infty$$

$$\leq 2\rho^{(k+m)/4} \sum_{i=-\infty}^{-\infty} \rho^{3|i|/4} \|\phi_i\|_\infty \sum_{j=-\infty}^{-\infty} \rho^{|j|/2} \|\phi_j\|_\infty.$$

The proof is complete. □

Thus $\{\Gamma_{k,m,x}(\theta)\}_{m\geq 0}$ is a uniform (w.r.t. $\theta \in G$) Cauchy sequence $\overline{\mathbb{P}}_{\theta^*}$-a.s. and in $L^1(\overline{\mathbb{P}}_{\theta^*})$, and $\{\Gamma_{k,m,x}(\theta)\}_{m\geq 0}$ converges as $m \to \infty$ uniformly w.r.t. $\theta$ $\overline{\mathbb{P}}_{\theta^*}$-a.s. and in $L^1(\overline{\mathbb{P}}_{\theta^*})$ to a random variable $\Gamma_{k,\infty}(\theta) \in L^1(\overline{\mathbb{P}}_{\theta^*})$ which does not depend on $x$ thanks to (44). By construction,

$$\text{var}_\theta \left[ \sum_{k=1}^n \phi(\theta, \mathbf{Z}_{k-q+1}^l) \Big| \overline{\mathbf{Y}}_0^n, X_0 = x_0 \right] = \sum_{k=1}^n \Gamma_{k,0,x_0}(\theta),$$

and the proof of Proposition 5 follows along the same lines as that of Proposition 4.

**Acknowledgment.** We are grateful to the anonymous referees who provided constructive criticism that improved the content and presentation of the paper.

R. Douc  
É. Moulines  
Ecole Nationale Supérieure  
 des Télécommunications  
Traitement du Signal et  
 des Images/CNRS URA 820  
46 Rue Barrault  
75634 Paris Cedex 13  
France  
e-mail: douc@tsi.enst.fr  
e-mail: moulines@tsi.enst.fr

T. Rydén  
Centre for Mathematical Sciences  
Lund University  
Box 118  
221 00 Lund  
Sweden  
e-mail: tobias@maths.lth.se